\documentclass[12pt]{article}
\usepackage{amsmath}
\usepackage[dvips]{graphicx}
\usepackage{textcomp}
\usepackage{amsbsy}
\usepackage{latexsym}
\usepackage[mathscr]{eucal}
\usepackage{amsfonts}
\newcommand{\Rbb}{\mathbb{R}}

\usepackage{amssymb}
\usepackage{algorithmic}
\usepackage{algorithm}
\usepackage{subfigure}

\usepackage[usenames,dvipsnames]{xcolor}
\usepackage{tikz,ifthen,fp}
\usepackage{pgfplots}
\usetikzlibrary{backgrounds,arrows}

\usepackage{hyperref}
\usepackage[makeroom]{cancel}

\usepackage{fullpage}

\newtheorem{lemma}{Lemma}[section]
\newtheorem{proposition}[lemma]{Proposition}

\newtheorem{remark}[lemma]{Remark}

\newtheorem{definition}[lemma]{Definition}

\def\beginproof{\noindent{\bf Proof:}~ }
\def\endproof{\hfill\rule{1.5mm}{1.5mm}\\[2mm]}

\begin{document}

\bibliographystyle{plain}

\title{Interpolation of inverse operators for preconditioning parameter-dependent equations
\thanks{This work was supported by the French National Research Agency (Grant ANR CHORUS MONU-0005)}
}
\author{ 
 Olivier ZAHM\footnotemark[2] ~and Anthony NOUY\footnotemark[2] \footnotemark[3]
}

\renewcommand{\thefootnote}{\fnsymbol{footnote}}
\footnotetext[2]{Ecole Centrale de Nantes, GeM, UMR CNRS 6183, France.}
\footnotetext[3]{Corresponding author (anthony.nouy@ec-nantes.fr).}

\renewcommand{\thefootnote}{\arabic{footnote}}

\maketitle
\date{}

\begin{abstract}
We propose a method for the construction of preconditioners of parameter-dependent matrices for the solution of large systems of parameter-dependent equations.
The proposed method is an interpolation of the matrix inverse based on a projection of the identity matrix with respect to the Frobenius norm. Approximations of the Frobenius norm using  random matrices  are introduced in order to handle large matrices. The resulting statistical estimators of the Frobenius norm yield quasi-optimal projections that are controlled with high probability. 
Strategies for the adaptive selection of interpolation points are then proposed for different objectives in the context of projection-based model order reduction methods: the improvement of residual-based error estimators, the improvement of the projection on a given reduced approximation space, or the re-use of computations for sampling based model order reduction methods. 
\end{abstract}

\section{Introduction}

This paper is concerned with the solution of large systems of parameter-dependent equations of the form
\begin{equation}\label{eq:Intro_AUB}
 A(\xi)u(\xi)=b(\xi),
\end{equation}
where $\xi$ takes values in some parameter set $\Xi$. Such problems occur in several contexts such as parametric analyses, optimization, control or uncertainty quantification, where  $\xi$ are random variables that parametrize model or data uncertainties.
The efficient solution of equation \eqref{eq:Intro_AUB} 
generally requires the construction of preconditioners for the operator $A(\xi)$, either for improving the performance of iterative solvers or for improving the quality of residual-based projection methods.

A basic preconditioner can be defined as the inverse (or any preconditioner) of the matrix $A(\bar{\xi})$ at some nominal parameter value $\bar{\xi}\in\Xi$ or as the inverse (or any preconditioner) of a mean value of $A(\xi)$ over $\Xi$ (see e.g. \cite{GHA96,ERN09}). When the operator only slightly varies over the parameter set $\Xi$, these parameter-independent preconditioners behave relatively well. However, for large variabilities, they are not able to provide a good preconditioning over the whole parameter set $\Xi$. A first attempt to construct a parameter-dependent preconditioner can be found in \cite{Desceliers2005}, where the authors compute through quadrature a polynomial expansion of the parameter-dependent factors of a LU factorization of $A(\xi)$. More recently, a linear Lagrangian interpolation of the matrix inverse has been proposed in  \cite{Chen14}. The generalization to any standard multivariate interpolation method is straightforward. However, standard approximation or interpolation methods require the evaluation of matrix inverses (or factorizations) for many instances of $\xi$ on a prescribed structured grid (quadrature or interpolation), that becomes prohibitive for large matrices and high dimensional parametric problems. 

In this paper, we propose an interpolation method for the inverse of matrix $A(\xi)$.
The interpolation is obtained by a projection of the inverse matrix on a linear span of samples of $A(\xi)^{-1}$ and takes the form 
\begin{equation*}
 P_m(\xi) = \sum_{i=1}^m \lambda_i(\xi) A(\xi_i)^{-1},
\end{equation*}
where $\xi_1,\hdots,\xi_m$ are $m$ arbitrary interpolation points in $\Xi$.  A natural interpolation could be obtained by minimizing  the condition number of $P_m(\xi)A(\xi)$ over the $\lambda_i(\xi)$, which is a Clarke regular strongly pseudoconvex optimization problem \cite{Ye2009}. However, the solution of this non standard optimization problem for many instances of $\xi$ is intractable and proposing an efficient solution method in a multi-query context remains a challenging issue. Here, the projection is defined as the minimizer of the Frobenius norm of $I-P_m(\xi)A(\xi)$, that is a quadratic optimization problem. Approximations of the Frobenius norm using random matrices are introduced in order to handle large matrices. These statistical estimations of the Frobenius norm  allow to obtain quasi-optimal projections that are controlled with high probability. Since we are interested in large matrices, $A(\xi_i)^{-1}$ are here considered as implicit matrices for which only efficient matrix-vector multiplications are available. Typically, a factorization (e.g. LU) of $A(\xi_i)$ is computed and stored. 
Note that when the storage of factorizations of several samples of the operator is unaffordable or when efficient preconditioners are readily available, one could similarly consider projections of the inverse operator on the linear span of preconditioners of samples of the operator. However, the resulting parameter-dependent preconditioner would be no more an interpolation of preconditioners. This straightforward extension of the proposed method is not analyzed in the present paper.  

The paper then presents several contributions in the context of projection-based model order reduction methods (e.g. Reduced Basis, Proper Orthogonal Decomposition (POD), Proper Generalized Decompositon) that rely on the projection of the solution $u(\xi)$ of \eqref{eq:Intro_AUB} on a low-dimensional approximation space. We first show how the proposed preconditioner can be used to define a Galerkin projection-based on the preconditioned residual, which can be interpreted as a Petrov-Galerkin projection of the solution with a parameter-dependent test space. 
Then, we propose adaptive construction of the preconditioner, based on an adaptive selection of interpolation points, for different objectives: (i) the improvement of error estimators based on preconditioned residuals, (ii) the improvement of the quality of projections on a given low-dimensional approximation space, or (iii) the re-use of computations for sample-based model order reduction methods. Starting from a $m$-point interpolation, these adaptive strategies consist in choosing a new interpolation point based on different criteria. In (i), the new point is selected for minimizing the distance between the identity and the preconditioned operator. In (ii), it is selected for improving the quasi-optimality constant of {Petrov-}Galerkin projections which measures how far the projection is from the best approximation on the reduced approximation space. In (iii), the new interpolation point is selected as a new sample determined for the approximation of the solution and not of the operator. The interest of the latter approach is that when direct solvers are used to solve equation \eqref{eq:Intro_AUB} at some sample points, the corresponding factorizations of the matrix can be stored and the preconditioner can be computed with a negligible additional cost. 

The paper is organized as follows. In Section \ref{sec:Frob_Interpolation} we present the method for the interpolation of the inverse of a parameter-dependent matrix. In Section \ref{sec:Model_reduction}, we show how the preconditioner can be used for the definition of a Petrov-Galerkin projection of the solution of \eqref{eq:Intro_AUB} on a given reduced approximation space, and we provide an analysis of the quasi-optimality constant of this projection. Then, different strategies for the selection of interpolation points for the preconditioner are proposed in Section \ref{sec:Greedy_Precond}. Finally, in Section   \ref{sec:Num_res}, numerical experiments will illustrate the efficiency of the proposed preconditioning strategies for different projection-based model order reduction methods. 

Note that the proposed preconditioner could be also used (a) for improving the quality of Galerkin projection methods where a projection of the solution $u(\xi)$ is searched on a subspace of functions of the parameters (e.g. polynomial or piecewise polynomial spaces) \cite{DEB01,MAT05,NOU09d}, or (b) for preconditioning iterative solvers for \eqref{eq:Intro_AUB}, in particular solvers based on low-rank truncations that require a low-rank structure of the preconditioner \cite{KHO11,MAT12,Giraldi:2014:tobeornottobe,Giraldi:2014}. These two potential applications are not considered here. 

\section{Interpolation of the inverse of a parameter-dependent matrix using Frobenius norm projection} \label{sec:Frob_Interpolation}

In this section, we propose a construction of an interpolation of the matrix-valued function $\xi\mapsto A(\xi)^{-1} \in\mathbb{R}^{n \times n}$ for given interpolation points 
$\xi_1,\hdots,\xi_m$ in $\Xi$.  We let $P_i= A(\xi_i)^{-1}$, $1\le i \le m$.
For large matrices, the explicit computation of $P_i$ is usually not affordable. 
Therefore, $P_i$ is here considered as an implicit matrix and we assume that the product of $P_i$ with a vector can be computed efficiently. In practice, factorizations of matrices $A(\xi_i)$ are stored. 
 
\subsection{Projection using Frobenius norm}\label{sec:Frob_proj}

We introduce the subspace $Y_m = \mathrm{span}\{P_1,\hdots,P_m\}$ of $\mathbb{R}^{n\times n}$. An approximation $P_m(\xi)$ of $A(\xi)^{-1}$ in $Y_m$ is then  defined by 
\begin{equation}\label{eq:FrobProjection}
 P_m(\xi) = \underset{P\in Y_m}{\mbox{argmin }} \| I -PA(\xi)  \|_F,
\end{equation}
where $I$ denotes the identity matrix of size $n$, and $\|\cdot\|_F$ is the Frobenius norm such that $\| B \|_F^2 = \langle B,B\rangle_F$ with $\langle B , C \rangle_F = \mbox{trace}(B^T C)$. Since $A(\xi_i)^{-1} \in Y_m$, we have the interpolation property $P_m(\xi_i) = A(\xi_i)^{-1}$, $1\le i\le m$. The minimization of $\| I - PA \|_F$ has been first proposed in  \cite{Grote1997} for the construction of a preconditioner $P$ in a subspace of matrices with given sparsity pattern (SPAI method). The following proposition gives some properties of the operator $P_m(\xi)A(\xi)$ (see Lemma 2.6 and Theorem 3.2 in \cite{Gonzalez2006}).

\begin{proposition}\label{prop:res_gonzalez}
  Let $P_m(\xi)$ be defined by \eqref{eq:FrobProjection}. We have
  \begin{equation}\label{eq:spectral_analysis}
  (1-\alpha_m(\xi))^2 \leq \| I-P_m(\xi)A(\xi) \|_F^2 \leq n(1-\alpha_m^2(\xi) ),
  \end{equation}
  where $ \alpha_m(\xi)$ is the lowest singular value of $P_m(\xi)A(\xi)$ verifying $0\le \alpha_m(\xi) \le 1$, with $P_m(\xi)A(\xi) = I$ if and only if $\alpha_m(\xi)=1$. Also, the following bound holds for the condition number of $P_m(\xi)A(\xi)$:
  \begin{equation}\label{eq:bound_kappa}
  \kappa(P_m(\xi)A(\xi)) \leq \frac{\sqrt{ n - (n-1) \alpha_m^2(\xi) }}{\alpha_m(\xi)}.
  \end{equation}
  Under the condition $\| I-P_m(\xi)A(\xi) \|_F < 1$, equations  
  \eqref{eq:spectral_analysis} and \eqref{eq:bound_kappa} imply that 
  \begin{equation*}
  \kappa(P_m(\xi)A(\xi)) \leq \frac{\sqrt{ n - (n-1) (1-\| I-P_m(\xi)A(\xi) \|_F)^2 }}{1 - \| I-P_m(\xi)A(\xi) \|_F}.
  \end{equation*}
\end{proposition}
For all $\lambda \in \mathbb{R}^m$, we have 
\begin{align*}
  \| I-\sum_{i=1}^m \lambda_i P_i A(\xi) \|_F^2 &= n - 2  \lambda ^T S(\xi) + \lambda^T  M(\xi)\lambda, 
\end{align*}
where the matrix $M(\xi) \in \mathbb{R}^{m\times m}$ and the vector $S(\xi)\in\mathbb{R}^m$ are given by
\begin{equation*}
 M_{i,j}(\xi) = \mbox{trace}( A^T(\xi) P_i^T P_j A(\xi) ) ~~~\mbox{ and }~~~
 S_{i}(\xi)   = \mbox{trace}( P_i A(\xi) ) .
\end{equation*}
Therefore, the solution of problem \eqref{eq:FrobProjection} is $P_m(\xi)=\sum_{i=1}^m \lambda_i(\xi)P_i$ with $\lambda(\xi)$ the solution of $M(\xi)\lambda(\xi) = S(\xi)$. When considering a small number $m$ of interpolation points, the computation time for solving this system of equations is negligible. However, the computation of $M(\xi)$ and $S(\xi)$ requires the evaluation of traces of matrices $A^T(\xi) P_i^T P_j A(\xi)$ and $P_i A(\xi)$ for all $1\leq i,j \leq m$. Since the $P_i$ are implicit matrices, the computation of such products of matrices is not affordable for large matrices. Of course, since $\mbox{trace}(B) = \sum_{i=1}^n e_i^TBe_i$, the trace of an implicit matrix $B$ could be obtained by computing the product of $B$ with the canonical vectors $e_1,\hdots,e_n$, but this approach is clearly not affordable for large $n$.

Hereafter, we propose an approximation of the above construction 
using an approximation of the Frobenius norm which requires less computational efforts.

\subsection{Projection using a Frobenius semi-norm} \label{sec:SemiFro}
Here, we define an approximation $P_m(\xi)$ of $A(\xi)^{-1}$ in $Y_m$ by
\begin{equation}\label{eq:SemiFrobProjection}
 P_m(\xi) = \underset{P\in Y_m}{ \mbox{argmin } } \| (I -PA(\xi))V  \|_F,
\end{equation}
where $V \in \Rbb^{n\times K}$, with $ K\le n$. $B\mapsto \| BV  \|_F$ defines a semi-norm on $\Rbb^{n\times n}$. Here, we assume that the linear map $P \mapsto PA(\xi)V$ is injective on $Y_m$ so that the solution of \eqref{eq:SemiFrobProjection} is unique. This requires $K\geq m$ and is satisfied when $\mathrm{rank}(V)\ge m$ and $Y_m$ is the linear span of linearly independent invertible matrices. Then, the solution $P_m(\xi) = \sum_{i=1}^m \lambda_i(\xi) P_i$ of \eqref{eq:SemiFrobProjection} is such that the vector $\lambda(\xi)\in\mathbb{R}^{m}$ satisfies $M^V(\xi)\lambda(\xi) = S^V(\xi)$, with 
\begin{equation}\label{eq:defMV_SV}
 M_{i,j}^V(\xi) = \mbox{trace}( V^T A^T(\xi) P_i^T P_j A(\xi) V ) ~~~\mbox{ and }~~~
 S_{i}^V(\xi)   = \mbox{trace}( V^TP_i A(\xi)V ) .
\end{equation}

The procedure for the computation of $M^V(\xi)$ and $S^V(\xi)$ is given in Algorithm \ref{alg:compute_MV_SV}. Note that only $mK$ matrix-vector products involving the implicit matrices $P_i$ are required. 

\begin{algorithm}[ht]
\begin{algorithmic}[1]
\REQUIRE $A(\xi),$ $ \{P_1,\hdots,P_m\}$ and $V = (v_{1},\hdots,v_{K})$
\ENSURE $M^V(\xi)$ and $S^V(\xi)$
\STATE Compute the vectors $w_{i,k} = P_i  A(\xi)v_k \in \Rbb^n$, for $1\le k \le K$ and $1\le i \le m$
\STATE Set $W_i = (w_{i,1},\hdots,w_{i,K}) \in \Rbb^{n\times K}$, $1\le i \le m$
\STATE Compute $M^V_{i,j}(\xi)=\mbox{trace}(W_i^T W_j)$ for $1\le i , j \le m$
\STATE Compute $S^V_{i}(\xi)=\mbox{trace}(V^T W_i)$ for $1\le i  \le m$
\end{algorithmic}
\caption{Computation of $M^V(\xi)$ and $S^V(\xi)$}
\label{alg:compute_MV_SV}
\end{algorithm}

Now the question is to choose a matrix $V$ such that $\| (I-PA(\xi))V \|_F$ provides a good approximation of $\| I-PA(\xi) \|_F$ for any $P\in Y_m$ and $\xi\in\Xi$.

\subsubsection{Hadamard matrices for the estimation of the Frobenius norm of an implicit matrix} \label{sec:Hadamard_trace}

Let $B$ an implicit $n$-by-$n$ matrix (consider $B=I-PA(\xi)$, with $P\in Y_m$ and $\xi\in\Xi$). Following \cite{Bekas2007}, we show how Hadamard matrices can be used for the estimation of the Frobenius norm of an implicit matrix. The goal is to find a matrix $V$ such that $\| BV \|_F$ is a good approximation of $\| B \|_F$. The relation $\| BV \|_F^2=\mbox{trace}(B^TBVV^T)$ suggests that $V$ should be such that $VV^T$ is as close as possible to the identity matrix. For example, we would like $V$ to minimize 
\begin{equation*}
 err(V)^2 =  \frac{1}{n(n-1)} \sum_{i=1}^n \sum_{j \neq i}^n (VV^T)^2_{i,j} = \frac{\| I-VV^T \|_F^2}{{n(n-1)}},
\end{equation*}
which is the mean square magnitude of the off-diagonal entries of $VV^T$. 
The bound $err(V)\geq \sqrt{(n-K)/((n-1)K)}$ is known to hold for any $V\in\mathbb{R}^{n\times K}$ whose rows have unit norm \cite{Welch1974}. Hadamard matrices can be used to construct matrices $V$ such that the corresponding error $err(V)$ is close to the bound, see \cite{Bekas2007}.

A Hadamard matrix $H_s$ is a $s$-by-$s$ matrix whose entries are $\pm1$, and which satisfies $H_sH_s^T=sI$ where $I$ is the identity matrix of size $s$. For example,
\begin{equation*}
 H_2 = \begin{pmatrix} 1&1 \\ 1&-1 \end{pmatrix} 
\end{equation*}
is a Hadamard matrix of size $s=2$. The Kronecker product (denoted by $\otimes$) of two Hadamard matrices is again a Hadamard matrix. Then it is possible to build a Hadamard matrix whose size $s$ is a power of 2 using a recursive procedure: $H_{2^{k+1}}=H_{2}\otimes H_{2^{k}}$. The $(i,j)$-entry of this matrix is $(-1)^{a^Tb}$, where $a$ and $b$ are the binary vectors such that $i=\sum_{k\geq0} 2^k a_k$ and $j=\sum_{k\geq0} 2^k b_k$. For a sufficiently large $s=2^k\geq \max(n,K)$, we define the \emph{rescaled partial Hadamard matrix} $V\in\mathbb{R}^{n\times K}$ as the first $n$ rows and the first $K$ columns of $H_s/\sqrt{K}$.

\subsubsection{Statistical estimation of the Frobenius norm of an implicit matrix}\label{sec:Stat_trace}

For the computation of the Frobenius norm of $B$, we can also use a statistical estimator as first proposed in \cite{Hutchinson1990}. The idea is to define a random matrix $V\in\mathbb{R}^{n\times K}$ with a suitable distribution law $\mathcal{D}$ such that $\| BV \|_F$ provides a controlled approximation of $\|B\|_F$ with high probability. 

\begin{definition}\label{def:JL_MATRIX}
 A distribution $\mathcal{D}$ over $\mathbb{R}^{n\times K}$ satisfies the $(\varepsilon,\delta)$-concentration property if for all $B\in\mathbb{R}^{n\times n}$,
 \begin{equation}\label{eq:JL_MATRIX}
  \mathbb{P}( | \| BV \|_F^2 -\| B \|_F^2| \geq \varepsilon\| B \|_F^2 ) \leq \delta.
 \end{equation}
\end{definition}

Two distributions $\mathcal{D}$ will be considered here.

(a) The \textit{rescaled Rademacher distribution}. Here the entries of $V\in\mathbb{R}^{n\times K}$ are independent and identically distributed with $V_{i,j} = \pm K^{-1/2}$ with probability $1/2$. According to Theorem 13 in \cite{Avron2010}, the rescaled Rademacher distribution satisfies the $(\varepsilon,\delta)$-concentration property for 
 \begin{equation}\label{eq:Kgeq_RAD}
  K\geq 6\varepsilon^{-2}\ln(2n/\delta).
 \end{equation}

(b) The \textit{subsampled Randomized Hadamard Transform distribution} (SRHT), first introduced in \cite{Ailon2006}. Here we assume that $n$ is a power of $2$. It is defined by $V = K^{-1/2} (RH_{n}D )^T \in\mathbb{R}^{n\times K}$ where
 \begin{itemize}
  \item $D\in\mathbb{R}^{n\times n}$ is a diagonal random matrix where $D_{i,i}$ are independent Rademacher random variables (i.e. $D_{i,i}=\pm 1$ with probability $1/2$),
  \item $H_n\in\mathbb{R}^{n\times n}$ is a Hadamard matrix of size $n$ (see Section \ref{sec:Hadamard_trace}),
  \item $R\in\mathbb{R}^{K\times n}$ is a subset of $K$ rows from the identity matrix of size $n$ chosen uniformly at random and without replacement.
 \end{itemize}
 In other words, we randomly select $K$ rows of $H_n$ without replacement, and we multiply the columns by $\pm K^{-1/2}$. We can find in \cite{Tropp2011,Boutsidis2013} an analysis of the SRHT matrix properties. In the case where $n$ is not a power of $2$, we define the partial SRHT (P-SRHT) matrix $V \in\mathbb{R}^{n\times K}$ as the first $n$ rows of a SRHT matrix of size $s\times K$, where  $s = 2^{\lceil \log_2(n)\rceil}$ is the smallest power of $2$ such that $n\le s<2n$. The following proposition shows that the (P-SRHT) distribution satisfies the $(\varepsilon,\delta)$-concentration property.

\begin{proposition}\label{prop:SRHT_conf_interval}
 The (P-SRHT) distribution satisfies the $(\varepsilon,\delta)$-concentration property for 
 \begin{equation}\label{eq:Kgeq_SRHT}
  K \geq 2(\varepsilon^2 - \varepsilon^3/3)^{-1}\ln(4/\delta)(1+\sqrt{8\ln(4n/\delta)})^2.
 \end{equation}
\end{proposition}

\beginproof
 Let $B\in\mathbb{R}^{n\times n}$. We define the square matrix $\widetilde B$ of size $s=2^{\lceil \log_2(n)\rceil}$, whose first $n\times n$ diagonal block is $B$, and $0$ elsewhere. Then we have $\|\widetilde B\|_F=\|B\|_F$. The rest of the proof is similar to the one of Lemma 4.10 in \cite{Boutsidis2013}.  We consider the events $A = \{ (1-\varepsilon) \| \widetilde B \|_F^2\leq \| \widetilde B V \|_F^2\leq (1+\varepsilon)\| \widetilde B \|_F^2 \}$ and $E = \{ \max_i \| \widetilde{B}DH_s^T e_i \|_2^2 \leq (1+\sqrt{8\ln(2s/\delta)})^2 \| \widetilde B \|_F^2 \}$, where $e_i$ is the $i$-th canonical vector of $\mathbb{R}^s$. The relation $\mathbb{P}(A^c) \leq \mathbb{P}(A^c | E) + \mathbb{P}(E^c)$ holds. Thanks to Lemma 4.6 in \cite{Boutsidis2013} (with $t=\sqrt{8\ln(2s/\delta)}$) we have $\mathbb{P}(E^c) \leq \delta/2$. Now, using the scalar Chernoff bound (Theorem 2.2 in \cite{Tropp2011} with $k=1$) we have 
 \begin{align*}
  \mathbb{P}(A^c|E) &= \mathbb{P}( \| \widetilde B V \|_F^2 \leq (1-\varepsilon) \| \widetilde B \|_F^2 |E) + \mathbb{P}( \| \widetilde B V \|_F^2 \geq (1+\varepsilon)\| \widetilde B \|_F^2|E) \\
  &\leq (e^{-\varepsilon} (1-\varepsilon)^{-1+\varepsilon}  )^{K(1+\sqrt{8\ln(2s/\delta)})^{-2}} + (e^{\varepsilon} (1+\varepsilon)^{-1-\varepsilon}  )^{K(1+\sqrt{8\ln(2s/\delta)})^{-2}}\\
  &\leq 2(e^{\varepsilon} (1+\varepsilon)^{-1-\varepsilon}  )^{K(1+\sqrt{8\ln(2s/\delta)})^{-2}} \leq 2 e^{K ( - \varepsilon^2/2 + \varepsilon^3/6) (1+\sqrt{8\ln(2s/\delta)})^{-2}}.
 \end{align*}
 The condition \eqref{eq:Kgeq_SRHT} implies $\mathbb{P}(A^c \vert E) \leq \delta/2$, and then $\mathbb{P}(A^c) \leq \delta/2 + \delta/2=\delta$, which ends the proof.
\endproof

Such statistical estimators are particularly interesting for that they provide approximations of the Frobenius norm of large matrices, with a number of columns $K$ for $V$ which scales as the logarithm of $n$, see \eqref{eq:Kgeq_RAD} and \eqref{eq:Kgeq_SRHT}. However, the concentration property \eqref{eq:JL_MATRIX} holds only for a given matrix $B$. The following proposition \ref{prop:concentration_SEV} extends these concentration results for any matrix $B$ in a given subspace. The proof is inspired from the one of Theorem 6 in \cite{Dasgupta09}. The essential ingredient is the existence of an $\varepsilon$-net for the unit ball of a finite dimensional space (see \cite{Bourgain89}). 

\begin{proposition}\label{prop:concentration_SEV}
 Let $V\in\mathbb{R}^{n\times K}$ be a random matrix whose distribution $\mathcal{D}$ satisfies the $(\varepsilon,\delta)$-concentration property, with $\varepsilon\le 1$. Then, for any $L$-dimensional subspace of matrices $M_L\subset \mathbb{R}^{n\times n} $ and for any $C>1$, we have
 \begin{equation}\label{eq:concentration_SEV}
  \mathbb{P}\big( | \| BV \|_F^2 -\| B \|_F^2| \geq \varepsilon(C+1)/(C-1) \| B \|_F^2 , \forall B\in M_L \big) \leq (9C/\varepsilon)^L \delta.
 \end{equation}
 
\end{proposition}

\beginproof
 We consider the unit ball $\mathcal{B}_L = \{ B\in M_L : \| B\|_F\leq 1 \}$ of the subspace $M_L$. It is shown in \cite{Bourgain89} that for any $\widetilde\varepsilon>0$, there exists a net $\mathcal{N}^{\widetilde\varepsilon}_L\subset \mathcal{B}_L$ of cardinality lower than $(3/\widetilde\varepsilon)^L$ such that 
 \begin{equation*}
  \min_{B_{\widetilde\varepsilon}\in \mathcal{N}^{\widetilde\varepsilon}_L} \| B-B_{\widetilde\varepsilon} \|_F \leq \widetilde\varepsilon,\quad \forall B\in \mathcal{B}_L.
 \end{equation*}
 In other words, any element of the unit ball $\mathcal{B}_L$ can be approximated by an element of $\mathcal{N}^{\widetilde\varepsilon}_L$ with an error less than $\widetilde\varepsilon$. Using the $(\varepsilon,\delta)$-concentration property and a union bound, we obtain
 \begin{equation}\label{eq:tmp2357}
  | \| B_{\widetilde\varepsilon}V \|_F^2 -\| B_{\widetilde\varepsilon} \|_F^2| \leq \varepsilon\| B_{\widetilde\varepsilon} \|_F^2 ,\quad \forall B_{\widetilde\varepsilon}\in \mathcal{N}^{\widetilde\varepsilon}_{L},
 \end{equation}
 with a probability at least $1-\delta (3/\widetilde\varepsilon)^L$. We now impose the relation $\widetilde\varepsilon = \varepsilon/(3C)$, where $C>1$. To prove \eqref{eq:concentration_SEV}, it remains to show that equation \eqref{eq:tmp2357}   implies \begin{equation}| \| BV \|_F^2 -\| B \|_F^2| \leq \varepsilon(C+1)/(C-1) \| B \|_F^2, \quad \forall B\in M_L.\label{eq:tmp2358}\end{equation}
 
 We define $B^* \in \arg\max_{B\in \mathcal{B}_L } | \| BV \|_F^2 -\| B \|_F^2|$. Let $B_{\widetilde\varepsilon}\in \mathcal{N}^{\widetilde\varepsilon}_L$ be such that $\| B^*-B_{\widetilde\varepsilon} \|_F \leq \widetilde\varepsilon$, and 
 $B^*_{\widetilde\varepsilon} = \arg\min_{B\in \text{span}(B_{\widetilde\varepsilon})} \| B^*-B \|_F$. Then we have $ \| B^*- B^*_{\widetilde\varepsilon} \|_F^2 = \| B^* \|_F^2-\| B^*_{\widetilde\varepsilon} \|_F^2\leq \widetilde\varepsilon^2$ and $\langle B^*-B^*_{\widetilde\varepsilon} ,  B^*_{\widetilde\varepsilon} \rangle=0$, where $\langle \cdot,\cdot \rangle$ is the inner product associated to the Frobenius norm $\|\cdot\|_F$. We have
 \begin{align*}
  \eta &:= | \| B^*V \|_F^2 -\| B^* \|_F^2| = | \| (B^*-B^*_{\widetilde\varepsilon})V + B^*_{\widetilde\varepsilon}V \|_F^2 -\| B^*-B^*_{\widetilde\varepsilon} + B^*_{\widetilde\varepsilon} \|_F^2| \\
  &= | \| (B^*-B^*_{\widetilde\varepsilon})V \|_F^2 +2\langle (B^*-B^*_{\widetilde\varepsilon})V ,  B^*_{\widetilde\varepsilon}V \rangle + \| B^*_{\widetilde\varepsilon}V  \|_F^2 -\| B^*-B^*_{\widetilde\varepsilon}  \|_F^2 - \| B^*_{\widetilde\varepsilon} \|_F^2 |  \\
  &\leq | \| (B^*-B^*_{\widetilde\varepsilon})V \|_F^2 - \| B^*-B^*_{\widetilde\varepsilon}  \|_F^2 | + | \| B^*_{\widetilde\varepsilon}V  \|_F^2  - \| B^*_{\widetilde\varepsilon} \|_F^2 | + 2 \| (B^*-B^*_{\widetilde\varepsilon})V \|_F \|  B^*_{\widetilde\varepsilon}V \|_F.
 \end{align*}
 We now have to bound the three terms in the previous expression. Firstly, since $( B^*-B^*_{\widetilde\varepsilon} )/\| B^*-B^*_{\widetilde\varepsilon} \|_F \in \mathcal{B}_L$, the relation $ | \| (B^*-B^*_{\widetilde\varepsilon})V \|_F^2 - \| B^*-B^*_{\widetilde\varepsilon}  \|_F^2 | \leq \| B^*-B^*_{\widetilde\varepsilon} \|_F ^2 \eta \leq \widetilde\varepsilon^2 \eta  $ holds. Secondly, \eqref{eq:tmp2357} gives $| \| B^*_{\widetilde\varepsilon}V  \|_F^2  - \| B^*_{\widetilde\varepsilon} \|_F^2 | \leq \varepsilon \| B^*_{\widetilde\varepsilon} \|_F^2 \leq \varepsilon $. Thirdly, by definition of $\eta$, we can write $ \| (B^*-B^*_{\widetilde\varepsilon})V \|_F^2 \leq (1+\eta) \| B^*-B^*_{\widetilde\varepsilon} \|_F^2 \leq \widetilde\varepsilon^2(1+\eta) $ and $ \|  B^*_{\widetilde\varepsilon}V \|_F^2 \leq (1+\varepsilon)\| B^*_{\widetilde\varepsilon} \|_F^2 \leq 1+\varepsilon $, so that we obtain $ 2 \| (B^*-B^*_{\widetilde\varepsilon})V \|_F \|  B^*_{\widetilde\varepsilon}V \|_F \leq 2\widetilde\varepsilon  \sqrt{1+\varepsilon}\sqrt{1+\eta}$. Finally, from \eqref{eq:tmp2357}, we obtain 
  \begin{equation}\label{eq:tmp9839}
  \eta \leq \widetilde\varepsilon^2 \eta + \varepsilon + 2\widetilde\varepsilon  \sqrt{1+\varepsilon}\sqrt{1+\eta}
 \end{equation}
 Since $\varepsilon\leq 1$, we have $\widetilde\varepsilon=\varepsilon/(3C)<1/3$. Then $\widetilde\varepsilon^2\leq \widetilde\varepsilon$ and $\sqrt{1+\varepsilon}\leq 3/2$, so that \eqref{eq:tmp9839} implies 
 \begin{align*}
  \eta &\leq \widetilde\varepsilon \eta + \varepsilon + 3\widetilde\varepsilon \sqrt{1+\eta}\leq \widetilde\varepsilon \eta + \varepsilon + 3\widetilde\varepsilon (1+\eta/2)  \leq 3 \widetilde\varepsilon \eta + \varepsilon + 3\widetilde\varepsilon,
 \end{align*}
 and then $\eta \leq (\varepsilon + 3\widetilde\varepsilon ) /(1-3\widetilde\varepsilon )%\leq \varepsilon(C+1)/(C-\varepsilon)
 \leq \varepsilon(C+1)/(C-1)$. By definition of $\eta$, we can write  
 $   | \| BV \|_F^2 -\| B \|_F^2|  \leq \varepsilon (C+1)/(C-1) $ for any $B\in \mathcal{B}_L$, that implies \eqref{eq:tmp2358}.
\endproof

\begin{proposition} \label{prop:QuasiOptimalite_seminorm}
 Let $\xi\in\Xi$, and let $P_m(\xi)\in Y_m$ be defined by \eqref{eq:SemiFrobProjection} where $V\in\mathbb{R}^{n\times K}$ is a realization of a rescaled Rademacher matrix with
 \begin{equation}\label{eq:Rademacher_K_SEV}
  K\geq 6\varepsilon^{-2} \ln( 2 n (9C/\varepsilon)^{m+1}/\delta) ,
 \end{equation}
 or a realization of a P-SRHT matrix with
 \begin{equation}\label{eq:SRHT_K_SEV}
 K \geq 2(\varepsilon^2 - \varepsilon^3/3)^{-1}\ln(4(9C/\varepsilon)^{m+1}/\delta)(1+\sqrt{8\ln(4n(9C/\varepsilon)^{m+1}/\delta)})^2
 \end{equation}
 for some $\delta>0$, $\varepsilon\leq 1$ and $C>1$. Assuming $\varepsilon'=\varepsilon(C+1)/(C-1)<1$, 
 \begin{equation}
  \| I-P_m(\xi)A(\xi) \|_F \leq \sqrt{\frac{1+\varepsilon'}{1-\varepsilon'}} \min_{P\in Y_m} \| I -PA(\xi)  \|_F
  \label{eq:quasi_optimality}
 \end{equation}
 holds with a probability higher than  $1-\delta$.
\end{proposition}

\beginproof
 Let us introduce the subspace $M_{m+1}=Y_mA(\xi)+\text{span}(I)$ of dimension less than $m+1$, such that $\{I-PA(\xi):P\in Y_m\} \subset M_{m+1}$. Then, we note that with the conditions \eqref{eq:Rademacher_K_SEV} or \eqref{eq:SRHT_K_SEV}, the distribution law $\mathcal{D}$ of the random matrix $V$ satisfies the $(\varepsilon,\delta (\varepsilon/(9C))^{m+1})$-concentration property. Thanks to Proposition \ref{prop:concentration_SEV}, the probability that
  $$ | \| (I-PA(\xi) )V \|_F^2 - \| I-PA(\xi) \|_F^2 | \leq \varepsilon' \| I-PA(\xi) \|_F^2 $$
 holds for any $P\in Y_m$ is higher than $1-\delta$. Then, by definition of $P_m(\xi)$ \eqref{eq:SemiFrobProjection}, we have with a probability at least $1-\delta$ that for any $P\in Y_m$, it holds 
 \begin{align*}
  \| I-P_m(\xi)A(\xi) \|_F &\leq \frac{1}{\sqrt{1-\varepsilon'}}  \| (I-P_m(\xi)A(\xi))V \|_F, \\
  &\leq \frac{1}{\sqrt{1-\varepsilon'}} \| (I-PA(\xi))V \|_F \leq \frac{\sqrt{1+\varepsilon'}}{\sqrt{1-\varepsilon'}} \| I-PA(\xi) \|_F.
 \end{align*} Then, taking the minimum over $P\in Y_m$, we obtain \eqref{eq:quasi_optimality}.
\endproof

Similarly to Proposition \ref{prop:res_gonzalez}, we obtain the following properties for $P_m(\xi)A(\xi)$, with $P_m(\xi)$ the solution of \eqref{eq:SemiFrobProjection}.

\begin{proposition}\label{prop:boundcond_seminorm}
Under the assumptions of Proposition \ref{prop:QuasiOptimalite_seminorm}, the inequalities
 \begin{align}\label{eq:spectral_analysis_V}
{ (1-\alpha_m(\xi))^2 }{(1-\varepsilon')^{-1}} &\leq \| (I-P_m(\xi)A(\xi) )V \|_F^2 \leq n\left( 1 - (1-\varepsilon')\alpha_m^2(\xi) \right)
 \end{align}
 and
\begin{align}   \kappa(P_m(\xi)A(\xi)) &\leq \alpha_m(\xi)^{-1} {\sqrt{ {n}{(1-\varepsilon')^{-1}} - (n-1) \alpha_m^2(\xi) }}{}\label{eq:bound_kappa_V}
 \end{align}
 hold with probability $1-\delta$, where $\alpha_m(\xi)$ is the lowest singular value of $P_m(\xi)A(\xi)$.
\end{proposition}

\beginproof
The optimality condition for $P_m(\xi)$ yields $\| (I-P_m(\xi)A(\xi) )V\|_F^2 = \| V \|_F^2 - \|P_m(\xi)A(\xi)V\|_F^2$. Since $P_m(\xi)A(\xi)\in M_{m+1}$ (where $M_{m+1}$ is the subspace introduced in the proof of Proposition \eqref{prop:QuasiOptimalite_seminorm}), we have 
\begin{equation}\label{eq:tmp5743}
 \| P_m(\xi)A(\xi)V\|_F^2 \geq (1-\varepsilon')\| P_m(\xi)A(\xi)\|_F^2
\end{equation}
with a probability higher than $1-\delta$. Using $\| V \|_F^2 = n$ (which is satisfies for any realization of the rescaled Rademacher or the P-SRHT distribution), we obtain $\| (I-P_m(\xi)A(\xi) )V\|_F^2 \le n - (1-\varepsilon') \|P_m(\xi)A(\xi)\|_F^2 $ with a probability higher than $1-\delta$. Then, $\| P_m(\xi)A(\xi) \|_F^2 \geq n \alpha_m(\xi)^2 $ yields the right inequality of \eqref{eq:spectral_analysis_V}. Following the proof of Lemma 2.6 in \cite{Gonzalez2006}, we have 
$(1-\alpha_m(\xi)^2) \leq \| I-P_m(\xi)A(\xi) \|_F^2$. Together with \eqref{eq:tmp5743}, it yields the left inequality of \eqref{eq:spectral_analysis_V}.
Furthermore, with probability $1-\delta$, we have $n-(1-\varepsilon') \| P_m(\xi)A(\xi) \|_F^2 \geq 0$. Since the square of the Frobenius norm of matrix $P_m(\xi)A(\xi) $ is the sum of squares of its singular values, we deduce 
$$
 (n-1) \alpha_m(\xi)^2 +  \beta_m(\xi)^2 \leq \| P_m(\xi)A(\xi) \|_F^2 \leq {n}{(1-\varepsilon')^{-1}}
$$
with a probability higher than $1-\delta$, where $\beta_m(\xi)$ is the largest singular value of $P_m(\xi)A(\xi)$. Then \eqref{eq:bound_kappa_V} follows from the definition of $\kappa(P_m(\xi)A(\xi) ) = {\beta_m(\xi)}/{\alpha_m(\xi)}$.
\endproof

\subsubsection{Comparison and comments}

\begin{figure}[h!]
   \centering
   \subfigure[$err(V)$ as function of $K$.]{\centering
     \includegraphics[width=.4\textwidth]{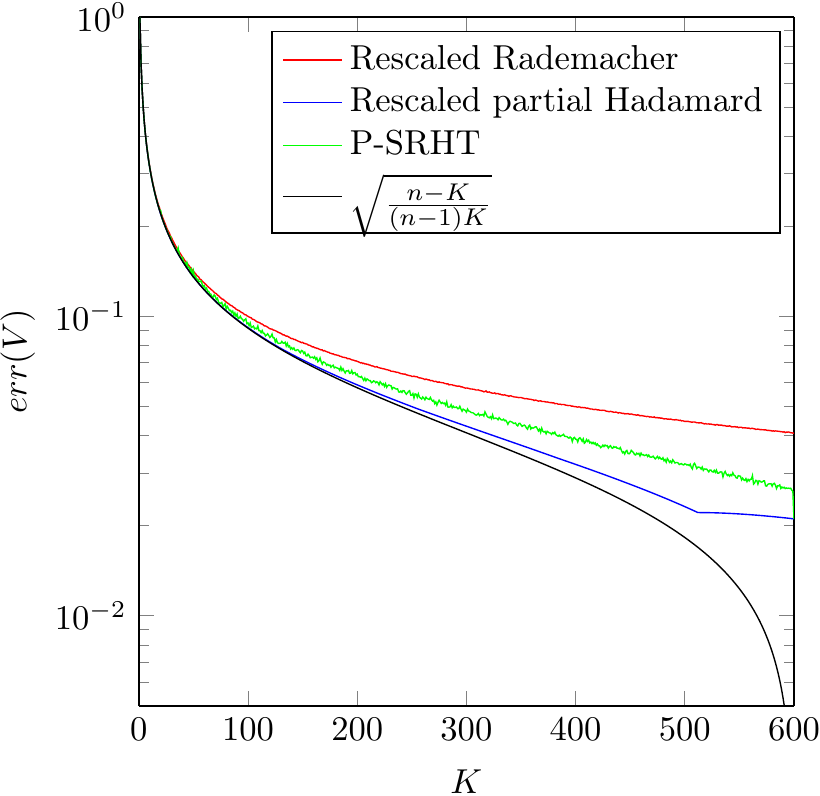}
     \label{fig:compare_Hadamard_MC}
   }~~~~~
   \subfigure[Distribution of the entries of $VV^T$ (in absolute value) where $V$ is the rescaled partial Hadamard matrix with $K=100$.]{\centering
     \includegraphics[width=.4\textwidth]{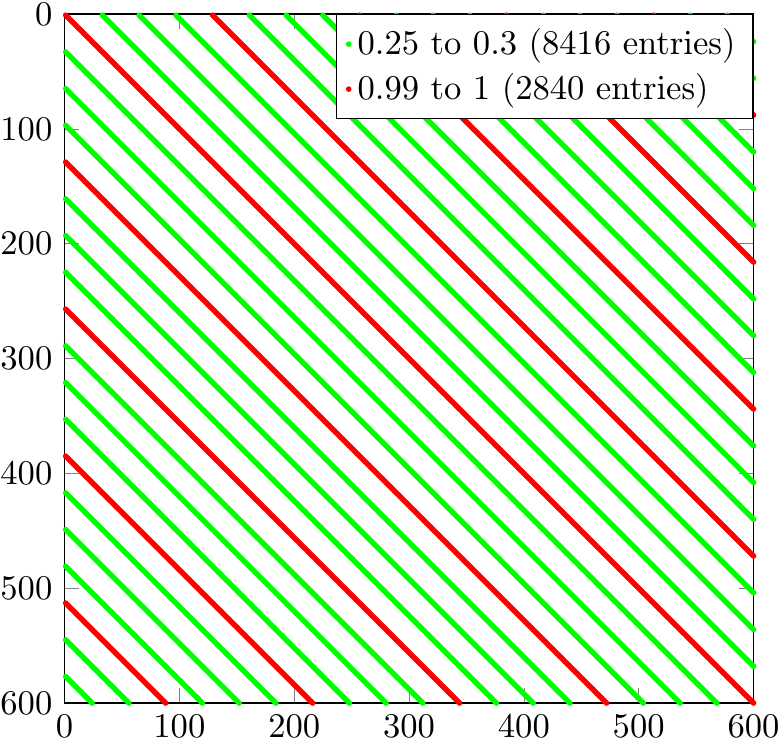}
     \label{fig:Hadamard_spy}
   }\\
   \subfigure[Distribution of the entries of $VV^T$ (in absolute value) where $V$ is a sample of the P-SRHT matrix with $K=100$.]{\centering
     \includegraphics[width=.4\textwidth]{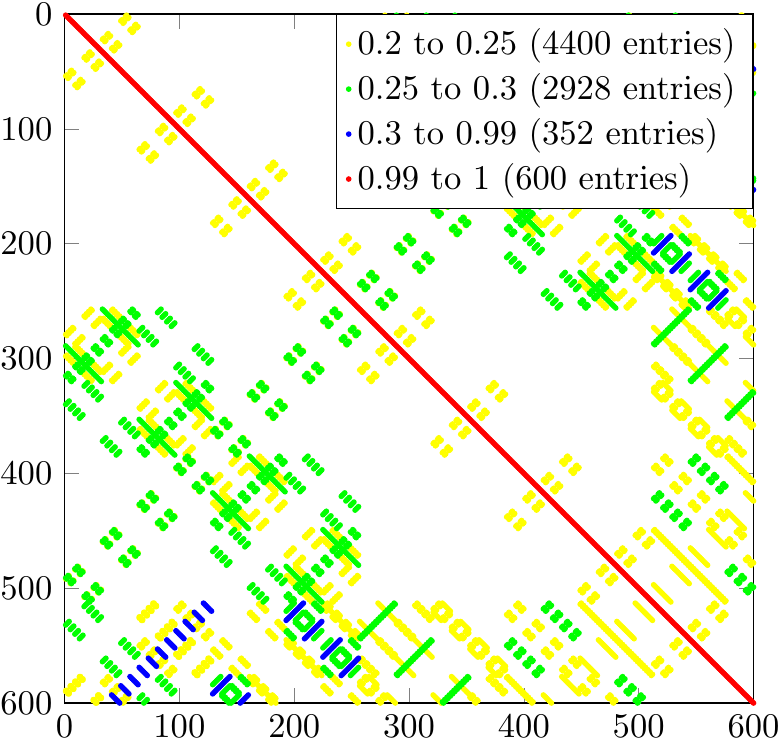}
     \label{fig:SRHT_spy}
     }~~~~~
   \subfigure[Distribution of the entries of $VV^T$ (in absolute value) where $V$ is a sample of the rescaled Rademacher matrix with $K=100$.]{\centering
     \includegraphics[width=.4\textwidth]{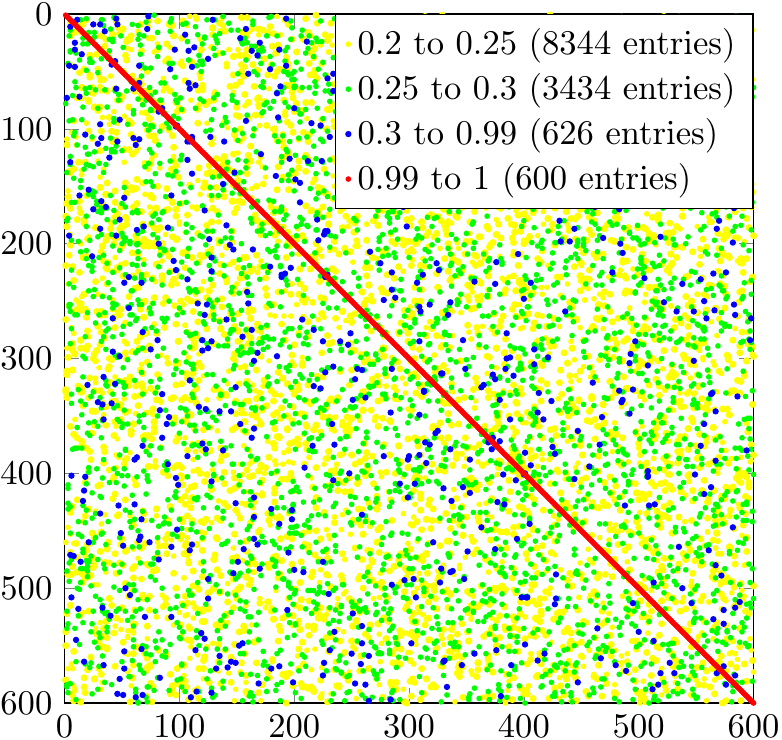}
     \label{fig:MC_spy}
     }
   \caption{Comparison between the rescaled partial Hadamard, the rescaled Rademacher and the P-SRHT matrix for the definition of matrix $V$, with $n=600$.}
   \label{fig:compare_Hadamard_MC_spy}
\end{figure}

We have presented different possibilities for the definition of $V$. The rescaled partial Hadamard matrices introduced in section \ref{sec:Hadamard_trace} have the advantage that the error $err(V)$ is close to the theoretical bound $\sqrt{(n-K)/((n-1)K)}$, see Figure \ref{fig:compare_Hadamard_MC} (note that the rows of $V$ have unit norm). Furthermore, an interesting property is that $VV^T$ has a structured pattern (see Figure \ref{fig:Hadamard_spy}). As noticed in \cite{Bekas2007}, when $K=2^q$ the matrix $VV^T$ have non-zero entries only on the $2^{qk}$-th upper and lower diagonals, with $k\geq0$. As a consequence, the error on the estimation of $\| B \|_F$ will be induced only by the non-zero off-diagonal entries of $B$ that occupy the $2^{qk}$-th upper and lower diagonals, with $k\geq1$. If the entries of $B$ vanish away from the diagonal, the Frobenius norm  is expected to be accurately estimated. Note that the P-SRHT matrices can be interpreted as a ``randomized version'' of the rescaled partial Hadamard matrices, and Figure \ref{fig:compare_Hadamard_MC} shows that the error $err(V)$ associated to the P-SRHT matrix behaves almost like the rescaled partial Hadamard matrix. Also, P-SRHT matrices yield a structured pattern for $VV^T$, see Figure \ref{fig:SRHT_spy}. The rescaled Rademacher matrices give higher errors $err(V)$ and yield matrices $VV^T$ with no specific patterns, see Figure \ref{fig:MC_spy}.

The advantage of using rescaled Rademacher matrices or P-SRHT matrices is that the quality of the resulting projection $P_m(\xi)$ can be controlled with high probability, provided that $V$ has a sufficiently large number of rows $K$ (see Proposition \ref{prop:QuasiOptimalite_seminorm}). Table \ref{tab:compare_RADE_SRHT} shows the theoretical value for $K$ in order to obtain the quasi-optimality result \eqref{eq:quasi_optimality} with $\sqrt{(1+\varepsilon')/(1-\varepsilon')}=10$ and $\delta=0.1\%$. It can be observed that $K$ grows very slowly with the matrix size $n$. Also, $K$ depends on $m$ linearly for the rescaled Rademacher matrices and quadratically for the P-SRHT matrices (see equations \eqref{eq:Rademacher_K_SEV} and \eqref{eq:SRHT_K_SEV}). However, these theoretical bounds for $K$ are very pessimistic, especially for the P-SRHT matrices. In practice, it can be observed that a very small value for $K$ may provide very good results (see Section \ref{sec:Num_res}). Also, it is worth mentioning that our numerical experiments do not reveal significant differences between the rescaled partial Hadamard, the rescaled Rademacher and the P-SRHT matrices.

\begin{table}[h!]
  \footnotesize
  \centering
  \subfigure[Rescaled Rademacher distribution.]{\centering
     \begin{tabular}{|l|c|c|c|c|c|} \cline{2-6}
      \multicolumn{1}{c|}{~}& $m=2$ & $m=5$ & $m=10$ & $m=20$ &  $m=50$ \\\hline
      $n=10^4$ & 239 & 363 & 567 & 972   & 2 185 \\\hline
      $n=10^6$ & 270 & 395 & 599 & 1 005 & 2 219 \\\hline
      $n=10^8$ & 301 & 427 & 632 & 1 038 & 2 253 \\\hline
     \end{tabular}
     \label{tab:compare_Rademacher}
   }

   \subfigure[P-SRHT distribution.]{\centering
     \begin{tabular}{|l|c|c|c|c|c|} \cline{2-6}
      \multicolumn{1}{c|}{~} & $m=2$ & $m=5$ & $m=10$ & $m=20$ &  $m=50$ \\\hline
      $n=10^4$ & 27 059 & 63 298 & 155 129 & 455 851 & 2 286 645 \\\hline
      $n=10^6$ & 30 597 & 69 129 & 164 750 & 473 011 & 2 326 301 \\\hline
      $n=10^8$ & 34 112 & 74 929 & 174 333 & 490 126 & 2 365 914 \\\hline
     \end{tabular}
     \label{tab:compare_SRHT}
   }
  \caption{Theoretical number of columns $K$ for the random matrix $V$ in order to ensure \eqref{eq:quasi_optimality}, with $\sqrt{(1+\varepsilon')/(1-\varepsilon')}=10$ and $\delta=0.1\%$. The constant $C$ has been chosen  in order to minimize $K$.}
  \label{tab:compare_RADE_SRHT}
\end{table}

\subsection{Ensuring the invertibility of the preconditioner for positive definite matrix}\label{sec:ensuring_invertibility}

Here, we propose a modification of the interpolation which ensures that $P_m(\xi)$  is invertible when 
$A(\xi)$ is positive definite.

Since $A(\xi_i)$ is positive definite, $P_i = A(\xi_i)^{-1}$ is positive definite. 
We introduce the vectors $\gamma^-\in\mathbb{R}^{m}$ and $\gamma^+\in\mathbb{R}^m$ whose components 
\begin{equation*}
 \gamma_i^- = \inf_{w\in\mathbb{R}^n} \frac{\langle P_i w,w \rangle}{\| w \|^2} >0 ~~\mbox{ and }~~\gamma_i^+ = \sup_{w\in\mathbb{R}^n} \frac{\langle P_i w,w \rangle}{\| w \|^2} <\infty
\end{equation*}
correspond respectively to the lowest and highest eigenvalues of the symmetric part of $P_i$. Then, for any 
$P=\sum_{i=1}^m \lambda_i P_i \in Y_m$, 
\begin{equation}\label{eq:precond_coercive}
 \inf_{w\in\mathbb{R}^n} \frac{\langle P w,w \rangle}{\| w \|^2} \geq \langle \lambda^+,\gamma^- \rangle - \langle \lambda^-,\gamma^+ \rangle,
\end{equation}
where $\lambda^+ \ge 0$ and $\lambda^- \ge 0$ are respectively the positive and negative parts of $\lambda = \lambda^+-\lambda^- \in \Rbb^m$.
As a consequence, if the right hand side of \eqref{eq:precond_coercive} is strictly positive, then $P$ is invertible. Furthermore, we have $ \|P\| \leq \langle \lambda^+ + \lambda^- , C \rangle $, where $C\in\mathbb{R}^m$ is the vector of component $C_i = \|P_i\|$, where $\| P_i\|$ denotes the operator norm of $P_i$. If we assume that $\langle \lambda^+,\gamma^- \rangle - \langle \lambda^-,\gamma^+ \rangle>0$, the condition number of $P$ satisfies 
\begin{equation*}
 \kappa(P) = \|P\|~\|P^{-1}\| \leq \|P\| \left( \inf_{w\in\mathbb{R}^n} \frac{\langle P w,w \rangle}{\| w \|^2} \right)^{-1} \leq \frac{\langle \lambda^+ + \lambda^- , C \rangle}{\langle \lambda^+,\gamma^- \rangle - \langle \lambda^-,\gamma^+ \rangle}.
\end{equation*}
It is then possible to bound $\kappa(P)$ by $\bar\kappa$ by imposing
\begin{equation*}
 \langle \lambda^+ + \lambda^- , C \rangle \leq \bar \kappa(\langle \lambda^+,\gamma^- \rangle - \langle \lambda^-,\gamma^+ \rangle),
\end{equation*}
which is a linear inequality constraint on $\lambda^+$ and $\lambda^-$. 
We introduce two convex subsets of $Y_m$ defined by 
\begin{subequations}
\begin{align}
 Y_m^{\bar \kappa} &=  \left\{ \sum_{i=1}^m \lambda_i^+ P_i - \sum_{i=1}^m \lambda_i^- P_i~:~
 \begin{matrix}
  \lambda_i^+\geq 0,~\lambda_i^-\geq 0\\
  \langle \lambda^+,\gamma^- \rangle - \langle \lambda^-,\gamma^+ \rangle \geq 0 \\
  \langle \lambda^+,\bar \kappa\gamma^- -C\rangle - \langle \lambda^-,\bar \kappa\gamma^+ +C \rangle \geq 0
 \end{matrix}
 \right\}, \nonumber \\
 Y_m^+ &=  \left\{ \sum_{i=1}^m \lambda_i P_i~:~  \lambda_i \geq 0\right\}. \nonumber 
\end{align}
\end{subequations}

From \eqref{eq:precond_coercive}, we have that any nonzero element of $Y_m^+$ is invertible, while any nonzero element of 
$Y_m^{\bar \kappa}$ is invertible and has a condition number lower than $\bar\kappa$. Under the condition  $\bar \kappa\geq\max_i {C_i}/{\gamma^-_i}$, we have
\begin{equation}\label{eq:YmP_YmK_Ym}
 Y_m^+ \subset Y_m^{\bar \kappa} \subset Y_m.
\end{equation}
Then definitions \eqref{eq:FrobProjection} and \eqref{eq:SemiFrobProjection} for the approximation $P_m(\xi)$ can be replaced respectively by 
\begin{subequations}
\begin{align}
 P_m(\xi) &= \underset{P\in Y_m^+ \mbox{ or }Y_m^{\bar \kappa}}{\mbox{argmin }} \| I -PA(\xi) \|_F ,\label{eq:FrobProjection_plusK}\\
 P_m(\xi) &= \underset{P\in Y_m^+ \mbox{ or }Y_m^{\bar \kappa}}{\mbox{argmin }} \| (I -PA(\xi) )V \|_F,\label{eq:FrobProjection_plusK_V}
\end{align}
\end{subequations}
which are quadratic optimization problems with linear inequality constraints. Furthermore, since $P_i\in Y_m^+$ for all $i$, all the resulting projections $P_m(\xi)$ 
interpolate $A(\xi)^{-1}$ at the points $\xi_1,\hdots,\xi_m$.

The following proposition shows that properties  \eqref{eq:spectral_analysis} and \eqref{eq:bound_kappa} still hold
for the preconditioned operator. 

\begin{proposition}
 The solution $P_m(\xi)$ of \eqref{eq:FrobProjection_plusK} is such that $P_m(\xi)A(\xi)$ satisfies \eqref{eq:spectral_analysis} and \eqref{eq:bound_kappa}. 
 Also, under the assumptions of Proposition \ref{prop:QuasiOptimalite_seminorm}, the solution $P_m(\xi)$ of \eqref{eq:FrobProjection_plusK_V} is such that $P_m(\xi)A(\xi)$ satisfies \eqref{eq:spectral_analysis_V} and \eqref{eq:bound_kappa_V} {with a probability higher than $1-\delta$}. 
\end{proposition}

\beginproof
Since $Y_m^+$ (or $Y_m^{\bar \kappa}$) is a closed and convex positive cone, the solution $P_m(\xi)$ of \eqref{eq:FrobProjection_plusK} is such that $trace((I-P_m(\xi)A(\xi))^T (P_m(\xi)-P)A(\xi)) \ge 0$ for all $P\in Y_m^+$ (or $Y_m^{\bar \kappa}$). Taking $P=2P_m(\xi)$ and $P=0$, we obtain that $\mbox{trace}((I-P_m(\xi)A(\xi))^T P_m(\xi)A(\xi))  =0$, which implies $ \| P_m(\xi)A(\xi) \|_F^2 = \mbox{trace}(P_m(\xi)A(\xi))$. We refer to the proof of Lemma 2.6 and  Theorem 3.2 in \cite{Gonzalez2006} to deduce \eqref{eq:spectral_analysis} and \eqref{eq:bound_kappa}. Using the same arguments, we prove that the solution $P_m(\xi)$ of \eqref{eq:FrobProjection_plusK_V} satisfies $\| P_m(\xi)A(\xi)V \|_F^2 = \mbox{trace}(V^T P_m(\xi)A(\xi)V)$, and then that \eqref{eq:spectral_analysis_V} and \eqref{eq:bound_kappa_V} hold {with a probability higher than $1-\delta$}.
\endproof

\subsection{Practical computation of the projection} \label{sec:practical_projection}

Here, we detail how to efficiently compute $M^V(\xi)$ and $S^V(\xi)$ given in equation \eqref{eq:defMV_SV} in a multi-query context, i.e. for several different values of $\xi$. The same methodology can be applied for computing $M(\xi)$ and $S(\xi)$. We assume that the operator $A(\xi)$ has an affine expansion of the form
\begin{align}\label{eq:A_affine}
  A(\xi) = \sum_{k=1}^{m_A} \Phi_k(\xi) A_k,
\end{align}
where the $A_k $ are matrices in $ \Rbb^{n\times n}$ and the $\Phi_k:\Xi\to \mathbb{R}$ are real-valued functions. Then $M^V(\xi)$ and $S^V(\xi)$ also have the affine expansions
\begin{subequations}\label{eq:MS_affine}
 \begin{align}
 M^V_{i,j}(\xi) &= \sum_{k=1}^{m_A}\sum_{l=1}^{m_A} \Phi_k(\xi)\Phi_l(\xi) ~\mbox{trace}( V^T A^T_k P_i^T P_j A_l V), \label{eq:MS_affine_a}\\
 S^V_{i}(\xi)   &= \sum_{k=1}^{m_A}\Phi_k(\xi) ~\mbox{trace}( V^T P_i A_k V) \label{eq:MS_affine_b},
 \end{align}
\end{subequations}
respectively. Computing the multiple terms of these expansions would require many computations of traces of implicit matrices and also, it would require the computation of the affine expansion of $A(\xi)$. Here, 
we use the methodology introduced in \cite{Casenave2014} for obtaining affine decompositions with a lower number of terms. These decompositions only require the knowledge of functions $\Phi_k$ in the affine decomposition \eqref{eq:A_affine}, and 
evaluations of $ M^V_{i,j}(\xi)$  and $ S^V_{i}(\xi)$ (that means evaluations of $A(\xi)$) at some selected points. 
We briefly recall this methodology. 

Suppose that $g:\Xi\rightarrow X$, with $X$ a vector space, has an affine decomposition $g(\xi) = \sum_{k=1}^m \zeta_k(\xi) g_k$, with $\zeta_k : \Xi\to \Rbb$ and $g_k\in X$.  
We first compute an interpolation of $\zeta(\xi) = (\zeta_1(\xi),\hdots,\zeta_m(\xi))$ under the form $\zeta(\xi) = \sum_{k=1}^{m_g} \Psi_k(\xi) \zeta(\xi_k^*)$, with $m_g\le m$, where $\xi_1^*, \hdots,\xi_{m_g}^*$ are interpolation points and $\Psi_1(\xi),\hdots,\Psi_{m_g}(\xi)$ the associated interpolation functions. Such an interpolation  can be computed with the Empirical Interpolation Method \cite{Maday2009} described in Algorithm \ref{alg:compute_EIM}. Then, we obtain 
an affine decomposition $g(\xi) = \sum_{k=1}^{m_g} \Psi_k(\xi) g(\xi^*_k)$ which can be computed from  
evaluations of $g$ at interpolation points $\xi_k^*$.
\begin{algorithm}[ht]
\begin{algorithmic}[1]
\REQUIRE $(\zeta_1(\cdot),\hdots,\zeta_m(\cdot))$
\ENSURE $\Psi_1(\cdot),\hdots, \Psi_k(\cdot)$ and $\xi^*_1,\hdots,\xi^*_k$
\STATE Define $R_1(i,\xi)=\zeta_i(\xi)$ for all $i,\xi$
\STATE Initialize $e=1$, $k=0$
\WHILE{$e \geq tolerance$ (in practice the machine precision)}
\STATE $k=k+1$
\STATE Find $(i^*_k,\xi^*_k) \in \underset{i,\xi}{\mbox{argmax}}~ \vert R_k(i,\xi) \vert$
\STATE Set the error to $e = \vert R_k(i^*_k,\xi^*_k) \vert$
\STATE Actualize $R_{k+1}(i,\xi) = R_k(i,\xi) - R_k(i,\xi^*_k)R_k(i^*_k,\xi)/R_k(i^*_k,\xi^*_k)$ for all $i,\xi$
\ENDWHILE
\STATE Fill in the $k$-by-$k$ matrix $Q$ : $Q_{i,j}=\zeta_{i^*_i}(\xi^*_j)$ for all $1\leq i,j \leq k$
\STATE Compute $\Psi_i(\xi) = \sum_{j=1}^k (Q^{-1})_{i,j} \zeta_{i^*_j}(\xi)$ for all $\xi$ and $1\leq i\leq k$
\end{algorithmic}
\caption{Empirical Interpolation Method (EIM).}
\label{alg:compute_EIM}
\end{algorithm}

Applying the above procedure to both $M^V(\xi)$ and $S^V(\xi)$, we obtain 
\begin{equation}
 M^V(\xi) \approx \sum_{k=1}^{m_M} \Psi_k(\xi) ~M^V(\xi^*_k),~~~
 S^V(\xi) \approx \sum_{k=1}^{m_S} \widetilde  \Psi_k(\xi) ~S^V(\widetilde \xi^*_k)  \label{eq:MS_interp_affine}.
\end{equation}

The first (so-called \emph{offline}) step consists in computing the interpolation functions $\Psi_k(\xi)$ and $ \widetilde \Psi_k(\xi)$ and associated  interpolation points $\xi^*_k$ and $ \widetilde \xi^*_k$ using Algorithm \ref{alg:compute_EIM} with input $\{\Phi_i\Phi_j\}_{1\leq i,j\leq m_A}$ and  $\{\Phi_i\}_{1\leq i\leq m_A}$ respectively, and then in computing matrices $M^V(\xi^*_k)$ and vectors $S^V(\widetilde \xi^*_k)$ using Algorithm \ref{alg:compute_MV_SV}.  The second (so-called \emph{online}) step simply consists in computing the matrix $M^V(\xi)$ and the vector $S^V(\xi)$ for a given value of $\xi$ using \eqref{eq:MS_interp_affine}.

\section{Preconditioners for projection-based model reduction}\label{sec:Model_reduction}

We consider a parameter-dependent  linear equation
\begin{equation}\label{eq:parametric_AUB}
 A(\xi) u(\xi) = b(\xi),
\end{equation}
with $A(\xi) \in \Rbb^{n\times n}$ and $b(\xi) = \mathbb{R}^n$. Projection-based model reduction consists in projecting the solution $u(\xi)$ onto a well chosen approximation space $X_r\subset X := \Rbb^n$ of low dimension $r\ll n$. Such projections are usually defined by imposing the residual of \eqref{eq:parametric_AUB} to be orthogonal to a so-called test space of dimension $r$. The quality of the projection on $X_r$ depends on the choice of the test space. The latter can be defined as the approximation space itself $X_r$, thus yielding the classical Galerkin projection. However when the operator $A(\xi)$ is ill-conditioned (for example when $A(\xi)$ corresponds to the discretization of non coercive or weakly coercive operators), this choice may lead to  projections that are far from optimal. Choosing the test space as $\{ R_X^{-1}A(\xi) v_r : v_r\in X_r \}$, where $R_X^{-1}A(\xi)$ is called the ``supremizer operator'' (see e.g. \cite{Rozza2007}), corresponds to a minimal residual approach, which may also results in projections that are far from optimal.  In this section, we show how the preconditioner $P_m(\xi)$ can be used for the definition of the test space. We also show how it can improve the quality of residual-based error estimates, which is a key ingredient for the construction of suitable approximation space $X_r$ in the context of the Reduced Basis method.

$X$ is endowed with the norm $\|\cdot\|_X$ defined by $\|\cdot\|_X^2 = \langle R_X\cdot,\cdot \rangle$, where $R_X$ is a symmetric positive definite matrix and $\langle \cdot,\cdot \rangle$ is the canonical inner product of $\mathbb{R}^n$. We also introduce the dual norm $\| \cdot \|_{X'}=\| R_X^{-1}\cdot \|_{X}$ such that for any $v,w\in X$ we have $|\langle v, w\rangle |\leq \|v\|_X~\|w\|_{X'}$.

\subsection{Projection of the solution on a given reduced subspace} \label{sec:MR_optimal_proj}

Here, we suppose that the approximation space $X_r$ has been computed by some model order reduction method.
The best approximation of $u(\xi)$ on $X_r$ is $u_r^*(\xi) = \arg\min_{v\in X_r}\| u(\xi)-v \|_X$ and is characterized by the orthogonality condition
\begin{equation}\label{eq:orthogonality_condition_proj_orth_u}
\langle u_r^{*}(\xi) - u(\xi), R_X v_r\rangle = 0,  \quad  \forall v_r\in X_r,
\end{equation}
or equivalently by the Petrov-Galerkin orthogonality condition
\begin{equation}\label{eq:orthogonality_condition_proj_orth_u_PG}
\langle A(\xi)u_r^{*}(\xi) - b(\xi), A^{-T}(\xi) R_X v_r\rangle = 0, \quad  \forall v_r\in X_r.
\end{equation}

Obviously the computation of test functions $A^{-T}(\xi)R_X v_r$ for basis functions $v_r$ of $X_r$ is prohibitive. By replacing $A(\xi)^{-1}$ by  $P_m(\xi)$, we obtain the feasible Petrov-Galerkin formulation 
\begin{equation}\label{eq:orthogonality_condition_proj_precond_u}
\langle A(\xi)u_r(\xi) - b(\xi), P_m^{T}(\xi)R_X v_r\rangle = 0, \quad  \forall v_r\in X_r.
\end{equation}

Denoting by $U \in \Rbb^{n\times r}$ a matrix whose range is $X_r$, the solution of \eqref{eq:orthogonality_condition_proj_precond_u} is $u_r(\xi)=Ua(\xi)$ where the vector $a(\xi)\in\mathbb{R}^r$ is the solution of
\begin{equation*}
 \big(U^T R_X P_m(\xi)A(\xi) U\big)a(\xi) = \big(U^T R_X P_m(\xi)b(\xi)\big).
\end{equation*} 

Note that \eqref{eq:orthogonality_condition_proj_precond_u} corresponds to the standard Galerkin projection when replacing $P_m(\xi)$ by $R_X^{-1}$. Indeed, the orthogonality condition \eqref{eq:orthogonality_condition_proj_precond_u} becomes $\langle A(\xi)u_r(\xi) - b(\xi), v_r\rangle = 0 $ for all $v_r\in X_r$.

\begin{remark}
 Here, the preconditioner $P_m(\xi)$ is used for the definition of the parameter-dependent test space $\{ P_m^{T}(\xi)R_X v_r : v_r \in X_r\}$ which defines the Petrov-Galerkin projection \eqref{eq:orthogonality_condition_proj_precond_u}. 
 However, $P_m(\xi)$ could also be used to construct preconditioners for the solution of the linear system $\big(U^T A(\xi) U\big)a(\xi) = \big( U^T b(\xi) \big)$ corresponding to the Galerkin projection on $X_r$. Following the idea proposed in \cite{ElmanPrecond}, such preconditoner can take the form $(U^T P_m(\xi)U)$, thus yielding the preconditioned reduced linear system
 $$
  \big(U^T P_m(\xi) U\big) \big(U^T A(\xi) U \big) a(\xi) = \big(U^T P_m(\xi) U\big) \big(U^T b(\xi) \big).
 $$
Such preconditioning strategy can be used to accelerate the solution of the reduced system of equations when using iterative methods. However, and contrarily to \eqref{eq:orthogonality_condition_proj_precond_u}, this strategy does not change the definition of $u_r(\xi)$, which is the standard Galerkin projection.
\end{remark}

We give now a quasi-optimality result for the approximation $u_r(\xi)$. This analysis relies on the notion of $\delta$-proximality introduced in \cite{Dahmen2011}.
\begin{proposition}\label{prop:delta_prox_PG}
 Let $\delta_{r,m}(\xi)\in [0,1]$ be defined by
 \begin{equation}\label{eq:delta}
  \delta_{r,m}(\xi) = \max_{ v_r\in X_r } \min_{w_r\in X_r} \frac{\| v_r-R_X^{-1}(P_m(\xi)A(\xi))^T R_X w_r \|_X}{\| v_r \|_X} .
 \end{equation}
 
 The solutions $u^*_r(\xi)\in X_r$ and $u_r(\xi)\in X_r$ of \eqref{eq:orthogonality_condition_proj_orth_u} and \eqref{eq:orthogonality_condition_proj_precond_u} satisfy 
 \begin{equation}\label{eq:control_PG_1}
  \| u_r^*(\xi) - u_r(\xi) \|_X \leq \delta_{r,m}(\xi) \| u(\xi) - u_r(\xi) \|_X.
 \end{equation}
 
 Moreover, if $\delta_{r,m}(\xi)<1$ holds, then 
 \begin{equation}\label{eq:control_PG_2}
  \| u(\xi) - u_r(\xi) \|_X \leq  {(1-\delta_{r,m}(\xi)^2)^{-1/2}} \| u(\xi) - u^*_r(\xi) \|_X.
 \end{equation}
\end{proposition}
\beginproof
The orthogonality condition \eqref{eq:orthogonality_condition_proj_orth_u} yields 
 \begin{equation*}
  \langle u_r^*(\xi) - u_r(\xi), R_X v_r \rangle = \langle u(\xi) - u_r(\xi), R_X v_r \rangle 
  = \langle b(\xi) - A(\xi)u_r(\xi),A^{-T}(\xi)R_X v_r \rangle
 \end{equation*}
 for all $v_r\in X_r$. Using  \eqref{eq:orthogonality_condition_proj_precond_u}, we have that for any $w_r\in X_r$,
 \begin{align*}
  \langle u_r^*(\xi) - u_r(\xi),R_X v_r \rangle &= \langle b(\xi) - A(\xi)u_r(\xi),A^{-T}(\xi)R_X v_r - P_m(\xi)^T R_X w_r\rangle ,\\
  &= \langle u(\xi) - u_r(\xi),R_X v_r - (P_m(\xi)A(\xi))^T R_X  w_r\rangle ,\\
  &\leq \| u(\xi) - u_r(\xi) \|_X~\|  R_X v_r - (P_m(\xi)A(\xi))^T R_X  w_r \|_{X'}\\
  &=    \| u(\xi) - u_r(\xi) \|_X~\|  v_r - R_X^{-1}(P_m(\xi)A(\xi))^T R_X  w_r \|_{X} .
 \end{align*}
 
 Taking the infimum over $w_r \in X_r$ and by the definition of $\delta_{r,m}(\xi)$, we obtain
 \begin{align*}
  \langle u_r^*(\xi) - u_r(\xi),R_X v_r \rangle &\leq \delta_{r,m}(\xi) \| u(\xi) - u_r(\xi) \|_X~\| v_r \|_{X}.
 \end{align*}
 
 Then, noting that $u_r^*(\xi) - u_r(\xi) \in X_r$, we obtain
   \begin{align*}
  \| u_r^*(\xi) - u_r(\xi) \|_X = \sup_{ v_r\in X_r} \frac{\langle u_r^*(\xi) - u_r(\xi),R_X v_r \rangle }{\| v_r \|_X}&\leq \delta_{r,m}(\xi) \| u(\xi) - u_r(\xi) \|_X,
 \end{align*}
 that is \eqref{eq:control_PG_1}. Finally, using orthogonality condition \eqref{eq:orthogonality_condition_proj_orth_u}, we have that 
 \begin{align*}
  \| u(\xi) - u_r(\xi) \|_X^2 &= \| u(\xi) - u^*_r(\xi) \|_X^2 + \| u^*_r(\xi) - u_r(\xi) \|_X^2 ,\\
  &\leq \| u(\xi) - u^*_r(\xi) \|_X^2 + \delta_{r,m}(\xi)^2 \| u(\xi) - u_r(\xi) \|_X^2,
 \end{align*}
from which we deduce \eqref{eq:control_PG_2} when  $\delta_{r,m}(\xi)<1$.
\endproof

An immediate consequence of Proposition \ref{prop:delta_prox_PG} is that when $\delta_{r,m}(\xi)=0$, the Petrov-Galerkin projection $u_r(\xi)$ coincides with the orthogonal projection $u^*_r(\xi)$. Following \cite{Dahmen2013}, we show in the following proposition that $\delta_{r,m}(\xi)$ can be computed by solving an eigenvalue problem of size $r$.
\begin{proposition}\label{prop:delta_computable}
We have $\delta_{r,m}(\xi) = \sqrt{1-\gamma}$, where $\gamma$ is the lowest eigenvalue of the generalized eigenvalue problem $ Cx=\gamma Dx $, with 
 \begin{align*}
  C &= U^T B(B^TR_X^{-1}B)^{-1}B^TU \in\mathbb{R}^{r\times r},\\
  D &= U^T R_X U \in\mathbb{R}^{r\times r},
 \end{align*}
where $B = (P_m(\xi)A(\xi))^T R_X U \in\mathbb{R}^{n\times r}$ and where $U \in  \Rbb^{n\times r}$ is a matrix whose range is $X_r$. 
\end{proposition}

\beginproof
 Since the range of $U$ is $X_r$, we have 
 \begin{equation*}
  \delta_{r,m}(\xi)^2= \max_{a\in \mathbb{R}^r} \min_{b\in \mathbb{R}^r} \frac{\| Ua-R_X^{-1}B b \|_X^2}{\| Ua \|_X^2}.
 \end{equation*}
 For any $a\in\mathbb{R}^r$, the minimizer $b^*$ of $\| U a-R_X^{-1}B b\|_X^2$ over $b\in \mathbb{R}^r$ is given by $b^* = (B^TR_X^{-1}B)^{-1}B^TUa$. Therefore, we have $\| Ua-R_X^{-1}Bb^* \|_{X}^2 = \| Ua \|_X^2 - \langle Ua , B b^* \rangle$, and 
 \begin{equation*}
  \delta_{r,m}^2(\xi)=  1-\inf_{a\in \mathbb{R}^r} \frac{\langle U^T B(B^TR_X^{-1}B)^{-1}B^TU a ,a \rangle}{\langle U^TR_X Ua,a\rangle},
 \end{equation*}
 which concludes the proof.
\endproof

\subsection{Greedy construction of the solution reduced subspace} \label{sec:MR_greedy_approx}

Following the idea of the Reduced Basis method \cite{Rozza2008,Veroy2003}, a sequence of nested approximation spaces $\{X_r\}_{r\ge 1}$ in $X$ can be constructed by a greedy algorithm such that $X_{r+1}=X_r + \mbox{span}(u(\xi_{r+1}^{RB}))$, where $\xi_{r+1}^{RB}$ is a point where the error of approximation of $u(\xi)$ in $X_r$ is maximal. An ideal greedy algorithm using the best approximation in $X_r$ and an exact evaluation of the projection error  is such that 
\begin{subequations}\label{eq:ideal_greedy}
\begin{align} 
 u_r^*(\xi) & \mbox{ is the orthogonal projection of $u(\xi)$ on $X_r$ defined by \eqref{eq:orthogonality_condition_proj_orth_u}}, \\ 
 \xi_{r+1}^{RB} &\in \underset{\xi\in\Xi}{\mbox{ argmax }} \| u(\xi) - u_r^*(\xi) \|_X. \label{eq:RBoptimal_xi}
\end{align}
\end{subequations}

This ideal greedy algorithm is not feasible in practice since $u(\xi)$ is not known. Therefore, we rather rely on a 
feasible weak greedy algorithm such that 
\begin{subequations}\label{eq:greedy-weak}
\begin{align}
 u_r(\xi) & \mbox{ is the Petrov-Galerkin projection of $u(\xi)$ on $X_r$ defined by \eqref{eq:orthogonality_condition_proj_precond_u}}, \\
 \xi_{r+1}^{RB} &\in \underset{\xi\in\Xi}{\mbox{ argmax }} \| P_m(\xi)( A(\xi)u_r(\xi)-b(\xi) )\|_X \label{eq:RBoptimal_xi_residual}.
\end{align}
\end{subequations}

Assume that 
\begin{equation*}
 \underline{\alpha}_m~\| u(\xi)-u_r(\xi) \|_X \leq \| P_m(\xi)(A(\xi)u_r(\xi)-b(\xi)) \|_X \leq \bar{\beta}_m ~\| u(\xi)-u_r(\xi) \|_X
\end{equation*}
holds with $\underline{\alpha}_m = \inf_{\xi\in\Xi} \alpha_m(\xi) >0$ and $\bar{\beta}_m   = \sup_{\xi\in\Xi} \beta_m(\xi) <\infty $, where $\alpha_m(\xi)$ and $\beta_m(\xi)$ are respectively the lowest and largest singular values of $P_m(\xi)A(\xi)$ with respect to the norm $\| \cdot\|_X$, respectively defined by the infimum and supremum of $\| P_m(\xi)A(\xi) v\|_X$ over $v \in X$ such that $\| v\|_X=1$.
Then, we easily prove that algorithm \eqref{eq:greedy-weak} is such that 
\begin{equation}
\| u(\xi^{RB}_{r+1}) - u_r(\xi^{RB}_{r+1}) \|_X \geq \gamma_m \max_{\xi\in \Xi} \| u(\xi) - u_r(\xi) \|_X,\label{eq:weak_greedy_selection}
\end{equation}
where $\gamma_m = \underline{\alpha}_m /~ \bar{\beta}_m \leq 1$ measures how far the selection of the new point is from the ideal greedy selection. Under condition \eqref{eq:weak_greedy_selection}, convergence results for this weak greedy algorithm can be found in  \cite{Binev2011,DeVore2013}.

We give now sharper bounds for the preconditoned residual norm that exploits the fact that the approximation $u_r(\xi)$ is the Petrov-Galerkin projection.

\begin{proposition}
 Let $u_r(\xi)$  be the Petrov-Galerkin projection of $u(\xi)$ on $X_r$ defined by \eqref{eq:orthogonality_condition_proj_orth_u_PG}. Then we have
 \begin{equation*}
  \alpha_{r,m}(\xi)~\| u(\xi)-u_r(\xi) \|_X \leq \| P_m(\xi)(A(\xi)u_r(\xi)-b(\xi)) \|_X \leq \beta_{r,m}(\xi) ~\| u(\xi)-u_r(\xi) \|_X ,
 \end{equation*}
 with
 \begin{subequations}
 \begin{align}
  \alpha_{r,m}(\xi) &= \inf_{v\in X}\sup_{w_r\in X_r} \frac{\| (P_m(\xi)A(\xi))^{T}R_X v \|_{X'}}{\| v - w_r \|_{X}} ,\nonumber\\ 
  \beta_{r,m}(\xi) &=\sup_{v\in X}\inf_{w_r\in X_r} \frac{\| (P_m(\xi)A(\xi))^T R_X (v-w_r) \|_{X'}}{\| v \|_X}. \nonumber
 \end{align}
 \end{subequations}
\end{proposition}

\beginproof
 For any $ v\in X$ and $ w_r\in X_r$ and according to \eqref{eq:orthogonality_condition_proj_precond_u}, we have 
 \begin{align*}
  \langle u(\xi) - u_r(\xi) ,R_X v \rangle &= \langle b(\xi)-A(\xi)u_r(\xi),A^{-T}(\xi)R_Xv - P_m^T(\xi) R_X w_r \rangle \\
  &= \langle P_m(\xi) (b(\xi)-A(\xi)u_r(\xi)),(P_m(\xi)A(\xi))^{-T}R_X v - R_X w_r \rangle \\
  &\leq \| R \|_{X} ~\|(P_m(\xi)A(\xi))^{-T}R_X v - R_X w_r \|_{X'},
  \end{align*}
  where $R(\xi):= P_m(\xi) (b(\xi)-A(\xi)u_r(\xi))$.
  Taking the infimum over $w_r\in X_r$, dividing by $\| v \|_X$ and taking the supremum over $ v\in X$, we obtain
  \begin{align*}
   \| u(\xi) - u_r(\xi) \|_X &\leq  \| R(\xi) \|_X \sup_{ v\in X}\inf_{w_r\in X_r} \frac{\| (P_m(\xi)A(\xi))^{-T}R_X v - R_X w_r\|_{X'}}{\| v \|_X}, \\
   &=  \| R(\xi) \|_X \sup_{v\in X}\inf_{w_r\in X_r} \frac{\| v - w_r\|_{X}}{\| (P_m(\xi)A(\xi))^{T}R_X v \|_{X'}}, \\
   &=  \| R(\xi) \|_X \left( \inf_{ v\in X}\sup_{w_r\in X_r} \frac{\| (P_m(\xi)A(\xi))^{T}R_X v \|_{X'}}{\| v - w_r \|_{X}} \right)^{-1},
  \end{align*}
  which proves the first inequality. Furthermore, for any $v\in X$ and $w_r \in X_r$, we have 
  \begin{align*}
   \langle P_m(\xi)(b(\xi) - A(\xi)u_r(\xi)) , R_X v \rangle &= \langle b(\xi) - A(\xi)u_r(\xi) ,P_m^T(\xi)R_X (v- w_r) \rangle \\
   &\leq \| u(\xi) - u_r(\xi) \|_X ~ \| (P_m(\xi)A(\xi))^TR_X (v-w_r) \|_{X'}.
  \end{align*}
  Taking the infimum over $w_r \in X_r$, dividing by $\|v\|_X$ and taking the supremum over $ v\in X$, we obtain the second inequality.
\endproof

Since $X_r\subset X_{r+1}$, we have $\alpha_{r+1,m}(\xi)\geq \alpha_{r,m}(\xi)\geq \alpha_{m}(\xi)$ and $\beta_{r+1,m}(\xi)\leq \beta_{r,m}(\xi)\leq \beta_{m}(\xi)$. Equation \eqref{eq:weak_greedy_selection} holds with $\gamma_m$ replaced by the 
parameter $\gamma_{r,m} = \underline{\alpha}_{r,m} /~ \bar{\beta}_{r,m}$. Since $\gamma_{r,m}$  increases with $r$, a reasonable expectation is that the convergence properties of the weak greedy algorithm will improve when $r$ increases.

\begin{remark}
 When replacing $P_m(\xi)$ by $R_X^{-1}$, the preconditioned residual norm $\| P_m(\xi)(A(\xi)u_r(\xi) -b(\xi)) \|_X$ turns out to be the residual norm $\| A(\xi)u_r(\xi) -b(\xi) \|_{X'}$, which is a standard choice in the Reduced Basis method for the greedy selection of points (with $R_X$ being associated with the natural norm on $X$ or with a norm associated with the operator at some nominal parameter value). This can be interpreted as a basic preconditioning method with a parameter-independent preconditioner.  
\end{remark}

\section{Selection of the interpolation points}\label{sec:Greedy_Precond}

In this section, we propose strategies for the adaptive selection of the interpolation points. For a given set of interpolation points $\xi_1,\hdots,\xi_m$, three different methods are proposed for the selection of a new interpolation point $\xi_{m+1}$. 
The first method aims at reducing uniformly the error between the inverse operator and its interpolation. The resulting interpolation of the inverse is pertinent for preconditioning iterative solvers or estimating errors based on preconditioned residuals. The second method aims at improving  Petrov-Galerkin projections of the solution of a parameter-dependent equation on a given approximation space. 
The third method aims at reducing the cost for the computation of the preconditioner by reusing operators computed when solving samples of a 
 parameter-dependent equation.

\subsection{Greedy approximation of the inverse of a parameter-dependent matrix}
A natural idea is to select a new interpolation point where the preconditioner $P_m(\xi)$ is not a good approximation of $A(\xi)^{-1}$.  Obviously, an ideal strategy for preconditioning would be to choose $\xi_{m+1}$ where the condition number of $P_m(\xi)A(\xi)$ is maximal. The computation of the condition number for many values of $\xi$ being computationaly expensive, one could use upper bounds of this condition number, e.g. computed using SCM \cite{Huynh2007c}.

Here, we propose the following selection method: given an approximation $P_m(\xi)$ associated with interpolation points $\xi_1,\hdots,\xi_m$, a new point $\xi_{m+1} $ is selected such that 
\begin{equation}\label{eq:greedy_precond_frob}
 \xi_{m+1} \in \underset{\xi\in\Xi}{\mbox{argmax}}~ \| (I - P_m(\xi)A(\xi))V \|_F,
\end{equation}
where the matrix $V$ is either the random rescaled Rademacher matrix, or the P-SRHT matrix (see Section \ref{sec:SemiFro}).
This adaptive selection of the interpolation points yields  the construction of an increasingsequence of subspaces $Y_{m+1} = Y_{m} + \mbox{span}(A(\xi_{m+1})^{-1})$ in $Y= \Rbb^{n\times n}$. This algorithm is detailed below. 
\begin{algorithm}[ht]
\begin{algorithmic}[1]
\REQUIRE $A(\xi), V, M$.
\ENSURE Interpolation points $\xi_1,\hdots,\xi_M$ and interpolation $P_M(\xi)$.
\STATE Initialize $P_0(\xi)=I$
\FOR{$m=0$ to $M-1$} 
\STATE Compute the new point  $\xi_{m+1}$ according to \eqref{eq:greedy_precond_frob}
\STATE Compute a factorization of $A(\xi_{m+1})$
\STATE Define $A(\xi_{m+1})^{-1}$ as an implicit operator 
\STATE Update the space $Y_{m+1}=Y_{m}+\mbox{span}(A(\xi_{m+1})^{-1})$
\STATE Compute  $P_{m+1}(\xi)= \arg\min_{P\in Y_{m+1}} \| (I-PA(\xi))V\|_F$
\ENDFOR
\end{algorithmic}
\caption{Greedy selection of interpolation points.}
\label{alg:compute_Greedy_Precond}
\end{algorithm}

 The following lemma interprets the above construction as a weak greedy algorithm.
\begin{lemma}\label{lem:weak-greedy}
 Assume that $A(\xi)$ satisfies $ \underline{\alpha}_0\|\cdot\| \leq\|A(\xi) \cdot \|\leq \bar{\beta}_0\|\cdot\|$ for all $\xi \in \Xi$, and let $P_m(\xi)$ be defined by \eqref{eq:SemiFrobProjection}. {Under the assumption that there exists $\varepsilon\in[0,1[$ such that
 \begin{equation}\label{eq:tmp9642}
  | \| (I-PA(\xi))V \|_F^2 - \| I-PA(\xi) \|_F^2 | \leq \epsilon \| I-PA(\xi) \|_F^2
 \end{equation}
 holds for all $\xi\in\Xi$ and $P\in Y_m$}, we have
 \begin{equation}
  \|  P_m(\xi_{m+1}) - A(\xi_{m+1})^{-1} \|_F \ge \gamma_{\varepsilon} \max_{\xi\in \Xi} \min_{P\in Y_m} \|  P - A(\xi)^{-1} \|_F,\label{eq:weak-greedy}
 \end{equation}
 with  $\gamma_{\varepsilon} = \underline{\alpha}_0\sqrt{1-\varepsilon}/(\bar{\beta}_0\sqrt{1+\varepsilon})$, and with $\xi_{m+1}$ defined by \eqref{eq:greedy_precond_frob}. 
\end{lemma}
\beginproof
 Since $\| BC \|_F \leq \| B \|_F \| C \|$ holds for any matrices $B$ and $C$, with $\| C \|$ the operator norm of $C$, we have for all $P\in Y$,
 \begin{align*}
   &\| A(\xi)^{-1} - P \|_F \leq \|  I - PA(\xi)   \|_F \| A(\xi)^{-1} \|\leq \underline{\alpha}_0^{-1} \|  I - PA(\xi)   \|_F , \\
   & \|  I - PA(\xi)   \|_F \leq \| A(\xi)^{-1} - P \|_F \| A(\xi) \|\leq \bar{\beta}_0 \|  A(\xi)^{-1} - P \|_F.
 \end{align*}
 Then, thanks to \eqref{eq:tmp9642} we have
 \begin{align*}
  & \| A(\xi)^{-1} - P \|_F \leq (\underline{\alpha}_0 \sqrt{1-\varepsilon})^{-1} \|  (I - PA(\xi))V   \|_F  \\
  \text{and }\quad  &  \| ( I - PA(\xi))V \|_F  \leq \bar{\beta}_0 \sqrt{1+\varepsilon} \|  A(\xi)^{-1} - P \|_F,
 \end{align*}
 which implies
 \begin{equation*}
  \frac{1}{\bar{\beta}_0\sqrt{1+\varepsilon}} \| (I -  P A(\xi) )V \|_F \leq  \|  A(\xi)^{-1} - P \|_F \leq  \frac{1}{\underline{\alpha}_0\sqrt{1-\varepsilon}} \| (I -  P A(\xi) )V \|_F.
 \end{equation*}
 We easily deduce that $\xi_{m+1}$ is such that \eqref{eq:weak-greedy} holds.
\endproof

\begin{remark}
The assumption \eqref{eq:tmp9642} of Lemma \ref{lem:weak-greedy} can be proved to hold with high probability in two cases. A first case is when $\Xi$ is a training set of finite cardinality, where the results of Proposition \ref{prop:QuasiOptimalite_seminorm} can be extended to any $\xi\in\Xi$ by using a union bound. We then obtain 
that  \eqref{eq:tmp9642} holds with a probability higher than $1-\delta (\#\Xi)$. A second case is when 
$A(\xi)$ admits an affine decomposition \eqref{eq:A_affine} with $m_A$ terms. Then the space $M_L = \text{span}\{I-PA(\xi) : \xi\in\Xi,  P \in Y_m\}$ is of dimension $L \le 1+m_A m$ and  Proposition \ref{prop:concentration_SEV} allows to prove that assumption  \eqref{eq:tmp9642} holds with high probability.
\end{remark}

The quality of the resulting spaces $Y_m$ have to be compared with the Kolmogorov $m$-width of the set $ A^{-1}(\Xi):=\{A(\xi)^{-1}:\xi\in\Xi \} \subset Y$, defined by
\begin{equation}\label{eq:kolmo_Am1}
 d_m (A^{-1}(\Xi) )_Y = \underset{\footnotesize \begin{matrix} Y_m\subset Y \\ \mbox{dim}(Y_m)=m \end{matrix}}{\min} ~ \underset{\xi\in\Xi}{\sup} ~ \underset{P\in Y_m}{\min}  \| A(\xi)^{-1} - P \|_F,
\end{equation}
which evaluates how well the elements of $A^{-1}(\Xi)$ can be approximated on a $m$-dimensional subspace of matrices. \eqref{eq:weak-greedy} implies that the following results holds (see  Corollary 3.3 in \cite{DeVore2013}):
\begin{equation*}
 \| A(\xi)^{-1}-P_m(\xi) \|_F = 
 \begin{cases}
 \mathcal{O}(m^{-a}) & \mbox{ if  } d_m ( A^{-1}(\Xi) )_Y = \mathcal{O}(m^{-a})\\ 
 \mathcal{O}(e^{-\tilde c m^{b}}) &\mbox{ if  } d_m ( A^{-1}(\Xi) )_Y = \mathcal{O}(e^{- c m^{b}})
 \end{cases},
\end{equation*}
 where $\tilde c > 0$ is a constant which depends on $c$ and $b$. That means that if the Kolmogorov $m$-width has an algebraic or exponential convergence rate, then the weak greedy algorithm yields an error $\| P_m(\xi)-A(\xi)^{-1} \|_F$ which has the same type of convergence. Therefore, the proposed interpolation method will present good convergence properties when  $d_m ( A^{-1}(\Xi) )_Y$ rapidly decreases with $m$.
 
\begin{remark}
 When the parameter set $\Xi $ is $ [-1,1]^d$ (or a product of compact intervals), an exponential decay can be obtained when $A(\xi)^{-1}$ admits an holomorphic extension to a domain in $\mathbb{C}^d$ containing $\Xi$ (see \cite{Chen14_2}).
\end{remark}
\begin{remark}
Note that here, there is no constraint on the minimization problem over $Y_m$ (either optimal subspaces or subspaces constructed by the greedy procedure), so that we have no guaranty that the resulting approximations $Y_m$ are invertible  (see Section \ref{sec:ensuring_invertibility}). 
\end{remark}
\subsection{Selection of points for improving the projection on a reduced space}\label{sec:greedy_delta}

We here suppose that we want to find an approximation of the solution $u(\xi)$ of a parameter-dependent equation \eqref{eq:parametric_AUB} onto 
a low-dimensional  approximation space $X_r$, using a Petrov-Galerkin orthogonality condition given by \eqref{eq:orthogonality_condition_proj_precond_u}. 
The best approximation is considered as the orthogonal projection defined by \eqref{eq:orthogonality_condition_proj_orth_u}. The quantity $\delta_{r,m}(\xi)$ defined by \eqref{eq:delta} controls the quality of the Petrov-Galerkin projection on $X_r$ (see Proposition \ref{prop:delta_prox_PG}). As indicated in Proposition \ref{prop:delta_computable}, $\delta_{r,m}(\xi)$ can be efficiently computed.  Thus, we propose the following selection strategy which aims at improving the quality of the Petrov-Galerkin projection: given a preconditioner $P_m(\xi)$ associated with interpolation points $\xi_1,\hdots,\xi_m$, the next point $\xi_{m+1}$ is  selected such that
\begin{equation}\label{eq:greedy_precond_delta}
  \xi_{m+1} \in \underset{\xi\in\Xi}{\mbox{argmax}}~ \delta_{r,m}(\xi).
\end{equation}

The resulting construction is described by Algorithm \ref{alg:compute_Greedy_Precond} with the above selection of $\xi_{m+1}$. 
Note that this strategy is closely related with \cite{Dahmen2013}, where the authors propose a greedy construction of a parameter-independent test space for Petrov-Galerkin projection, with a selection of basis functions based on an error indicator similar to $\delta_{r,m}(\xi)$.

\subsection{Re-use of factorizations of operator's evaluations - Application to reduced basis method}\label{sec:recycling}
When using a sample-based approach for solving a parameter-dependent equation \eqref{eq:parametric_AUB}, the linear system is solved for many values of the parameter $\xi$. When using a direct solver for solving a linear system for a given $\xi$, a factorization of the operator is usually available and can be used for improving a preconditioner for the solution of subsequent linear systems. 

We here describe this idea in the particular context of greedy algorithms for Reduced Basis method, where the interpolation points $\xi_1,\hdots,\xi_r$ for the interpolation of the inverse $A(\xi)^{-1}$ are taken as the evaluation points $\xi_1^{RB},\hdots,\xi_{r}^{RB}$  for the solution. At iteration $r$, having a preconditioner $P_r(\xi)$ and an approximation $u_r(\xi)$, a new interpolation point is defined such that  
\begin{equation*}
 \xi_{r+1}^{RB} \in \underset{\xi\in\Xi}{\mbox{ argmax }} \| P_r(\xi)( A(\xi)u_r(\xi)-b(\xi) )\|_X.
\end{equation*}

Algorithm \ref{alg:greedy_simultaneous} describes this strategy.

\begin{algorithm}[ht]
\begin{algorithmic}[1]
\REQUIRE $A(\xi),b(\xi)$, $V$, and $R$.
\ENSURE Approximation $u_{R}(\xi)$.
\STATE Initialize $u_0(\xi)=0$, $P_{0}(\xi)=I$
\FOR{$r=0$ to $R-1$}
\STATE Find  $\xi^{RB}_{r+1} \in \arg\max_{\xi\in \Xi} \| P_{r}(\xi)(A(\xi)u_r(\xi) - b(\xi)) \|_X $
\STATE Compute a factorization of $A(\xi^{RB}_{r+1})$
\STATE Solve the linear system $v_{r+1}=A(\xi^{RB}_{r+1})^{-1} b(\xi^{RB}_{r+1})$
\STATE Update the approximation subspace $X_{r+1} = X_{r} + \mbox{span}(v_{r+1})$
\STATE Define the implicit operator $P_{r+1} = A(\xi^{RB}_{r+1})^{-1}$
\STATE Update the space $Y_{r+1}$ (or $Y_{r+1}^+$)
\STATE Compute the preconditioner : $P_{r+1}(\xi) = \underset{P\in Y_{r+1} (\mbox{or }Y_{r+1}^+)}{\mbox{argmin}} \| (I-PA(\xi))V \|_F$
\STATE Compute the Petrov-Galerkin approximation $u_{r+1}(\xi)$ of $u(\xi)$ on  $X_{r+1}$ using equation \eqref{eq:orthogonality_condition_proj_precond_u}
\ENDFOR
\end{algorithmic}
\caption{Reduced Basis method with re-use of operator's factorizations.}
\label{alg:greedy_simultaneous}
\end{algorithm}

\section{Numerical results} \label{sec:Num_res}

\subsection{Illustration on a one parameter-dependent model}

\begin{figure}[h!]
   \centering
   \subfigure[$f$]{\centering
     \includegraphics[width=.23\textwidth]{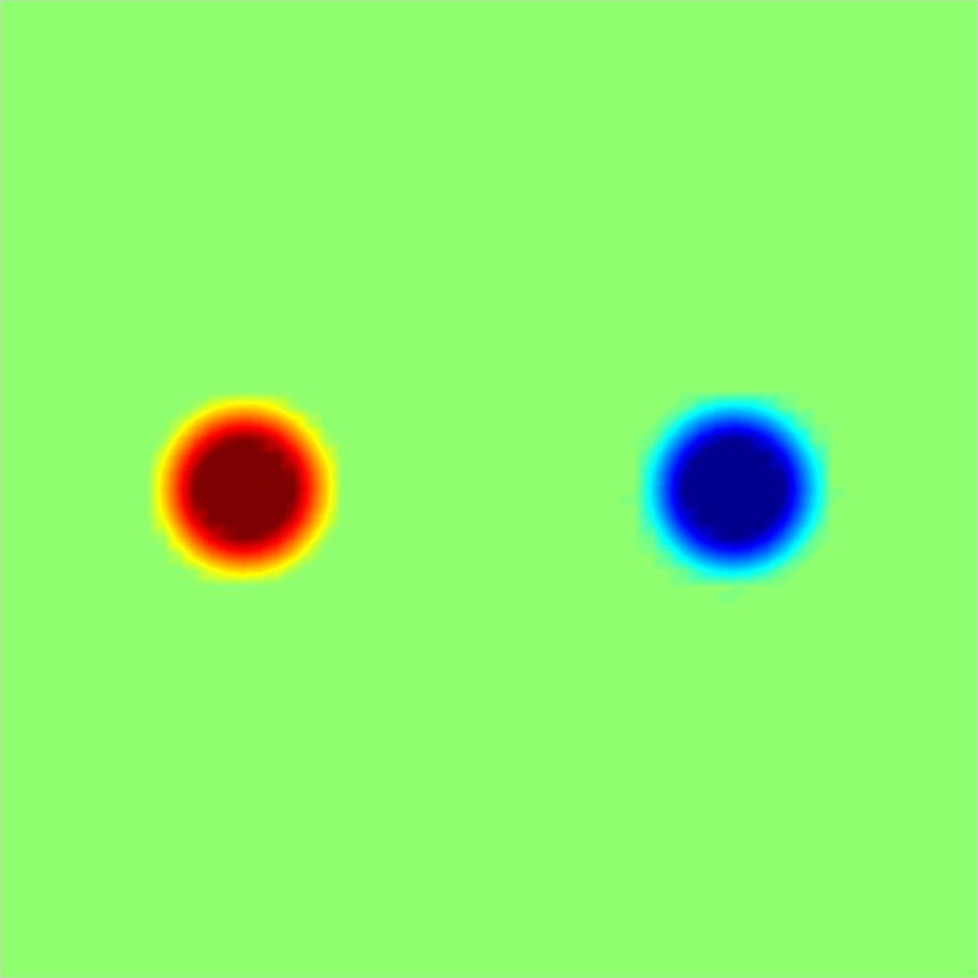} \label{fig:add_rot_source}
     }
   \subfigure[$u(0.05)$]{\centering
     \includegraphics[width=.23\textwidth]{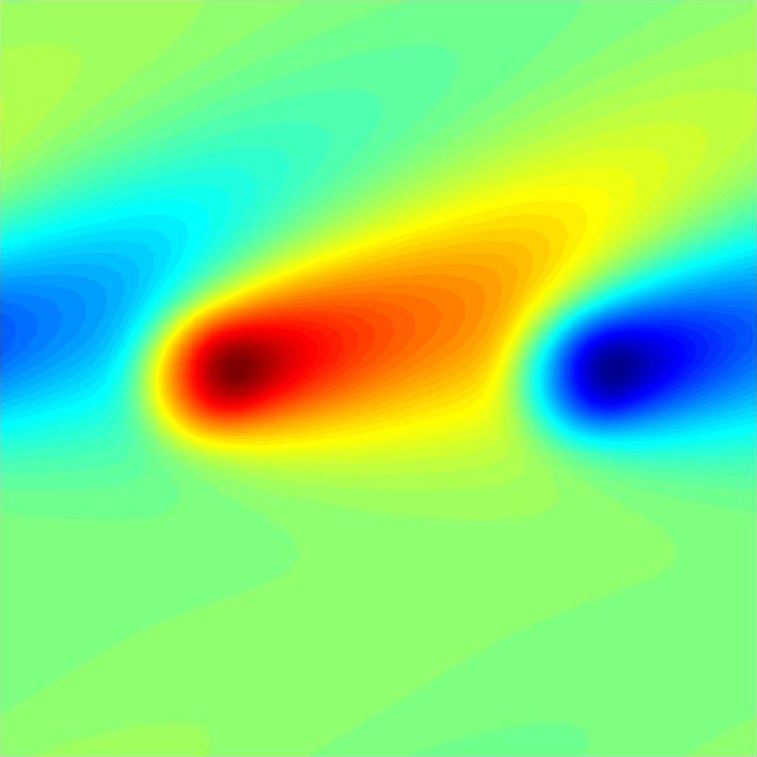}\label{fig:add_rot_sample1}
     }
   \subfigure[$u(0.2)$]{\centering
     \includegraphics[width=.23\textwidth]{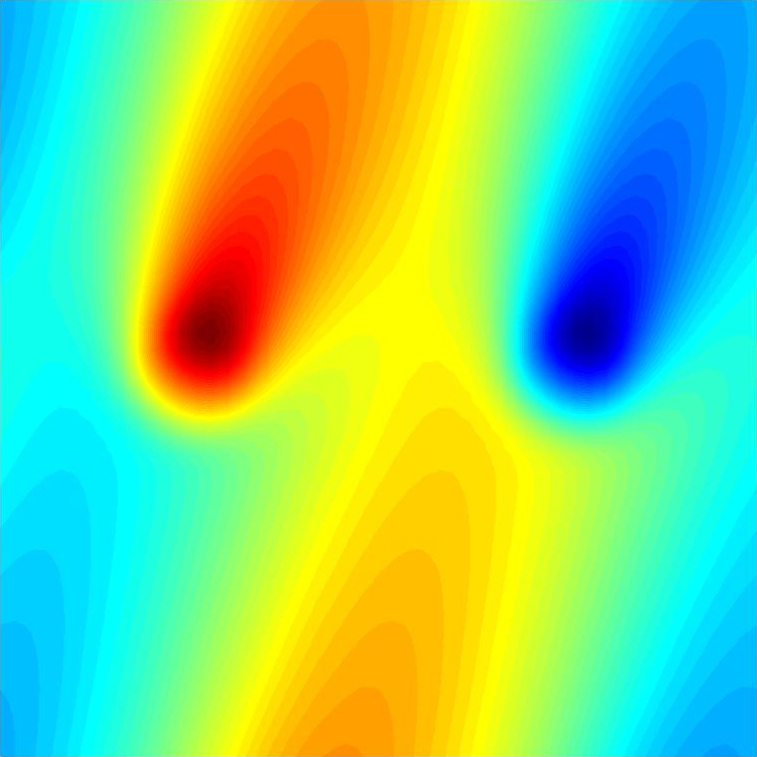}\label{fig:add_rot_sample2}
     }
   \subfigure[$u(0.8)$]{\centering
     \includegraphics[width=.23\textwidth]{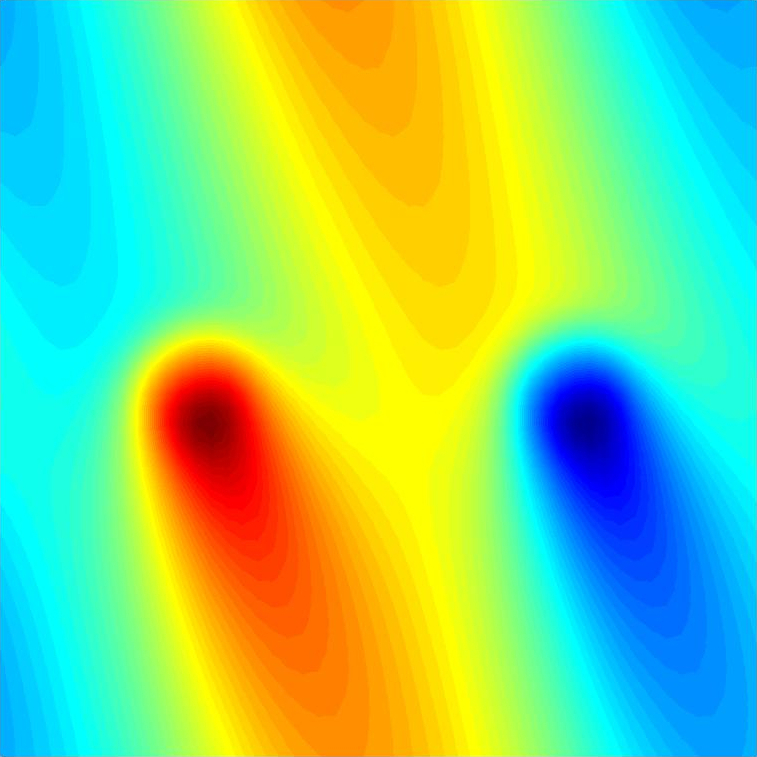}\label{fig:add_rot_sample3}
     }
   \caption{Plot of the source term $f$ (a) and 3 samples of the solution corresponding to parameter values $\xi=0.05$ (b), $\xi=0.2$ (c) and $\xi=0.8$ (d) respectively.}
   \label{fig:add_rot}
\end{figure}

In this section we compare the different interpolation methods on the following one parameter-dependent advection-diffusion-reaction equation:
\begin{equation}
 -\Delta u + v(\xi) \cdot \nabla u +u= f 
\end{equation}
defined over a square domain $\Omega = [0,1]^2$ with periodic boundary conditions. The advection vector field $v(\xi)$ is spatially constant and depends on the parameter $\xi$ that takes values in $[0,1]$:
$v(\xi)=D\cos(2\pi\xi)e_1 + D\sin(2\pi\xi) e_2 $, with $D=50$ and $(e_1,e_2)$ the canonical basis of $\mathbb{R}^2$. $\Xi$ denotes a uniform grid of 250 points on $[0,1]$. The source term $f$  is represented in Figure \ref{fig:add_rot_source}. We introduce a finite element approximation space of dimension $n=1600$ with piecewise linear approximations on a regular mesh of $\Omega$. 
The mesh P\'eclet number takes moderate values (lower than one), so that a standard Galerkin projection without stabilization is here sufficient. The Galerkin projection yields the linear system of equations $A(\xi)u(\xi)=b$, with 
\begin{align*}
 A(\xi)=A_0 + \cos(2\pi\xi) A_1 + \sin(2\pi\xi) A_2,
\end{align*}
where the matrices $A_0$, $A_1$, $A_2$ and the vector $b$ are given by 
\begin{align*}
 (A_0)_{i,j} &= \int_{\Omega} \nabla \phi_i \cdot\nabla \phi_j + \phi_i\phi_j  ~,~
 (A_1)_{i,j} = \int_{\Omega}  \phi_i (e_1 \cdot \nabla \phi_j )  \\
 (A_2)_{i,j} &= \int_{\Omega}  \phi_i (e_2 \cdot \nabla \phi_j )   ~,~
 (b)_i= \int_\Omega \phi_i f ,
\end{align*}
where $\{\phi_i\}_{i=1}^n$ is the basis of the finite element space. Figures \ref{fig:add_rot_sample1}, \ref{fig:add_rot_sample2} and \ref{fig:add_rot_sample3} show three samples of the solution.

\subsubsection{Comparison of the interpolation strategies}

We first choose arbitrarily 3 interpolation points ($\xi_1=0.05$, $\xi_2=0.2$ and $\xi_3=0.8$) and show the benefits of using the Frobenius norm projection for the definition of the  preconditioner. For the comparison, we consider the Shepard and the nearest neighbor interpolation strategies. Let $\|\cdot\|_\Xi$ denote a norm on the parameter set $\Xi$. The Shepard interpolation method is an inverse weighted distance  interpolation:
\begin{equation*}
 \lambda_i(\xi) = 
\left\{
\begin{array}{cl}
  \displaystyle \frac{\| \xi - \xi_i \|^{-s}_\Xi}{\sum_{j=1}^{m}\| \xi - \xi_j \|^{-s}_\Xi} &\mbox{ if } \xi\neq\xi_i \\
  1 & \mbox{ if } \xi = \xi_i 
\end{array}
\right. ,
\end{equation*}
where $s>0$ is a parameter. Here we take  $s=2$. The nearest neighbor interpolation method consists in choosing the value taken by the nearest interpolation point, that means $\lambda_i(\xi)=1$ for some $ i \in \arg\min_j \| \xi-\xi_j \|_\Xi$, and $\lambda_j(\xi)=0$ for all $ j\neq i$.

Concerning the Frobenius norm projection on $Y_m$ (or $Y_m^+$), we first construct the affine decomposition of $M(\xi)$ and $S(\xi)$ as explained in Section \ref{sec:practical_projection}. The interpolation points $ \xi^*_k$ (resp. $\widetilde \xi^*_k$) given by the EIM procedure for $M(\xi)$ (resp. $S(\xi)$) are  $\{0.0;0.25;0.37;0.56;0.80\}$ (resp. $\{0.0;0.25;0.62\}$). The number of terms $m_M=5$ in the resulting affine decomposition of $M(\xi)$ (see equation \eqref{eq:MS_interp_affine}) is less than the expected number $m_A^2=9$ (see equation \eqref{eq:MS_affine_a}). Considering the functions $\Phi_1(\xi)=1$, $\Phi_2(\xi)=\cos(2\pi\xi)$, $\Phi_3(\xi)=\sin(2\pi\xi)$, and thanks to relation $\cos^2=1-\sin^2$, the space
\begin{equation*}
 \mbox{span}_{i,j}\{\Phi_i\Phi_j\} = \mbox{span}\{1 ,\cos,\sin,\cos\sin,\cos^2,\sin^2\} = \mbox{span}\{1 ,\cos,\sin,\cos\sin,\cos^2\}
\end{equation*}
is of dimension $m_M=5$. The EIM procedure automatically detects the redundancy in the set of functions and reduces the number of terms in the decomposition \eqref{eq:MS_interp_affine}. Then, since the dimension $n$ of the discretization space is reasonable, we compute the matrices $M(\xi^*_k)$ and the vectors $S(\widetilde \xi^*_k)$ using equation \eqref{eq:MS_interp_affine}.

\begin{figure}[h!]
   \centering
   \subfigure[Nearest neighbor]{\centering
     \includegraphics[width=.28\textwidth]{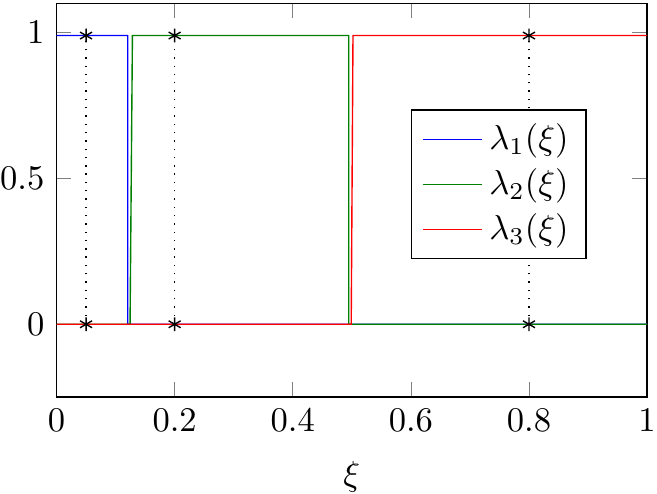}
     }
   \subfigure[Shepard (with $s=2$)]{\centering
     \includegraphics[width=.28\textwidth]{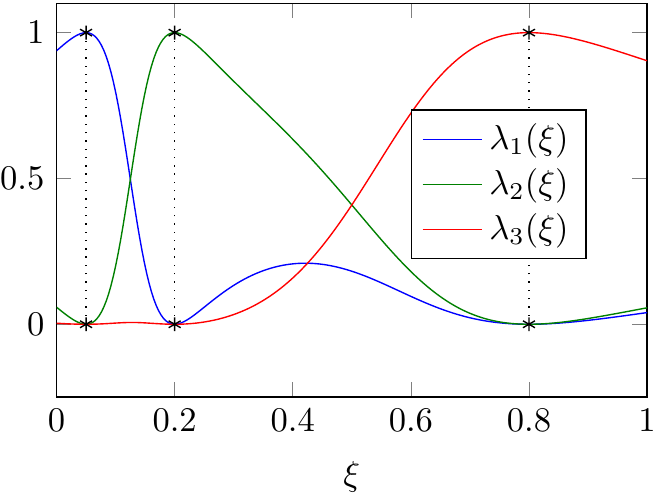}
     }
     \\
   \subfigure[Projection on $Y_m$]{\label{fig:add_rot_interpolation_function_c}\centering
     \includegraphics[width=.28\textwidth]{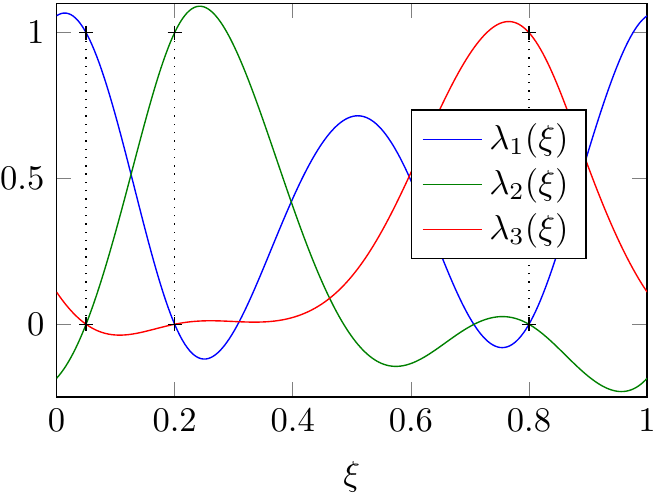}
     }
   \subfigure[Projection on $Y_m^+$]{\label{fig:add_rot_interpolation_function_d}\centering
     \includegraphics[width=.28\textwidth]{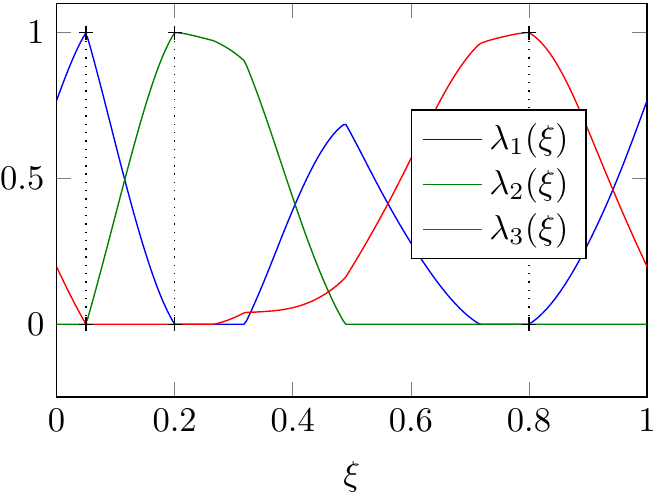}
     }
   \subfigure[Projection on $Y_m^{\bar \kappa}$]{\label{fig:add_rot_interpolation_function_e}\centering
     \includegraphics[width=.28\textwidth]{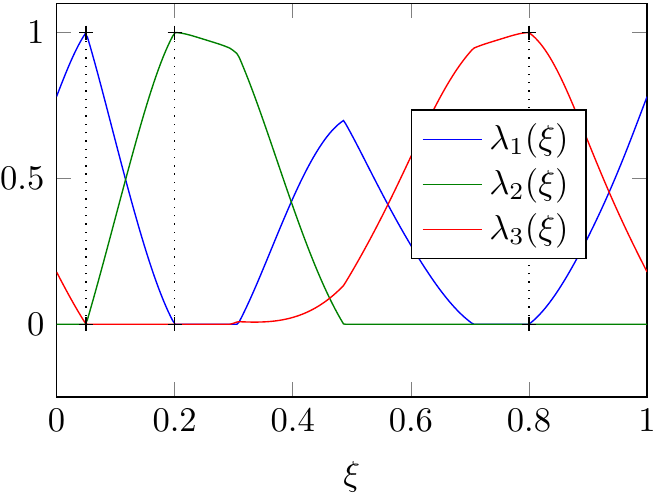}
     }
   \caption{Interpolation functions $\lambda_i(\xi)$ for different interpolation methods.}
   \label{fig:add_rot_interpolation_function}
\end{figure}

The functions $\lambda_i(\xi)$ are plotted on Figure \ref{fig:add_rot_interpolation_function} for the proposed interpolation strategies. It is important to note that contrary to the Shepard or the nearest neighbor method, the Frobenius norm projection (on $Y_m$ or $Y_m^+$) leads to periodic interpolation functions, \textit{i.e.} $\lambda_i(\xi=0)=\lambda_i(\xi=1)$. This is consistent with the fact that the application $\xi \mapsto A(\xi)$ is $1$-periodic. The Frobenius norm projection automatically detects such a feature.

\begin{figure}[h!]
  \centering
  \includegraphics[width=0.85\textwidth]{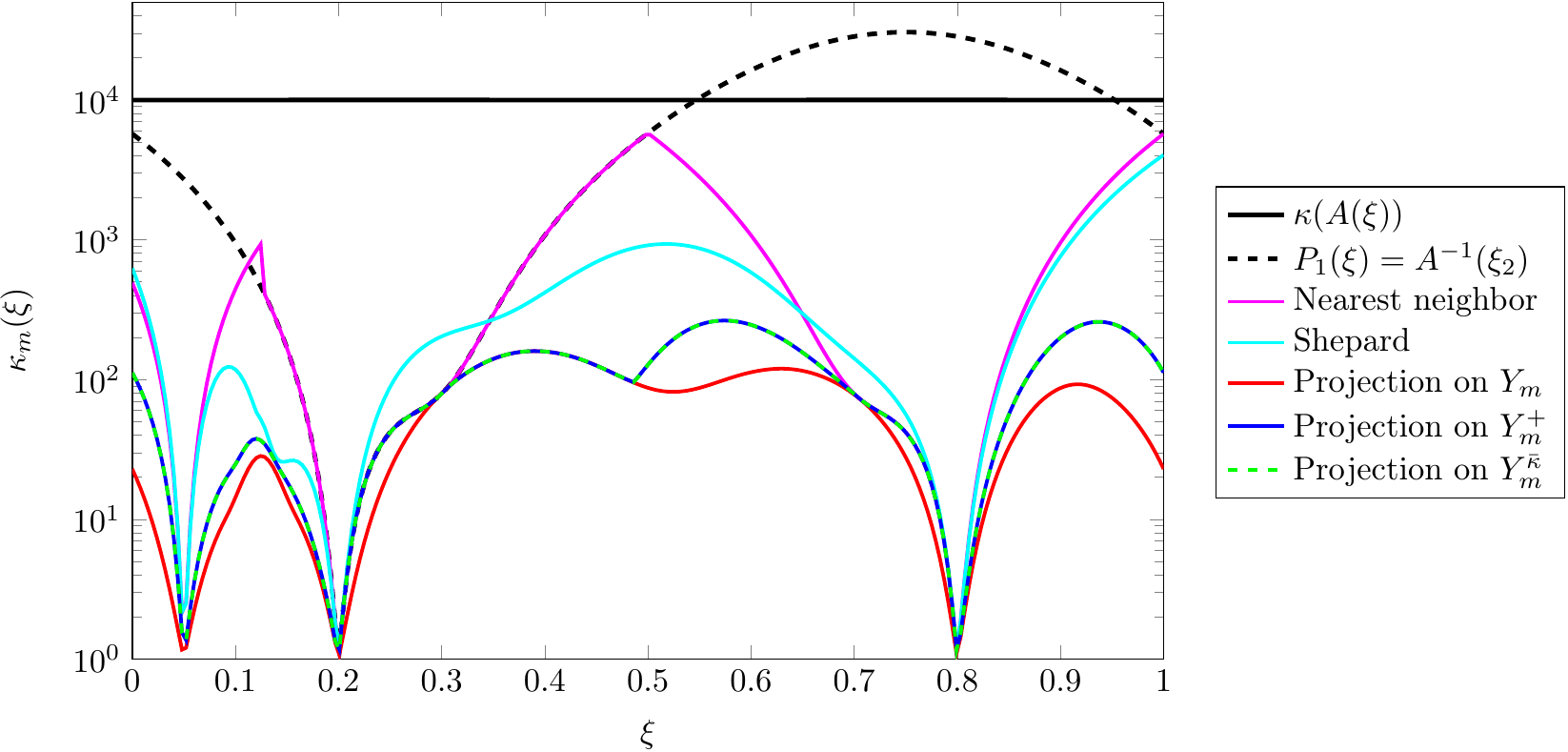}
  \caption{Condition number of $P_m(\xi)A(\xi)$ for different interpolation strategies. The condition number of $A(\xi)$ is given as a reference.}
  \label{fig:add_rot_condition_number}
\end{figure}

Figure \ref{fig:add_rot_condition_number} shows the condition number $\kappa_m(\xi)$ of $P_m(\xi)A(\xi)$ with respect to $\xi$. We first note that for the constant preconditioner $P_1(\xi)=A(\xi_2)^{-1}$, the resulting condition number is higher than the one of the non preconditioned matrix $A(\xi)$ for $\xi\in[0.55;0.95]$. We also note that the interpolation strategies based on the Frobenius norm projection lead to better preconditioners than the Shepard and nearest neighbor interpolation strategies. When considering the projection on $Y_m^+$ and $Y_m^{\bar \kappa}$ (with $\bar \kappa=5\times10^{4}$ such that \eqref{eq:YmP_YmK_Ym} holds), the resulting condition number is roughly
 the same, so as the interpolation functions of Figures \ref{fig:add_rot_interpolation_function_d} and \ref{fig:add_rot_interpolation_function_e}. Since the projection on $Y_m^{\bar \kappa}$ requires the expensive computation of the constants $\gamma^+$, $\gamma^-$ and $C$ (see Section \ref{sec:ensuring_invertibility}), we prefer to simply use the projection on $Y_m^+$ in order to ensure the preconditioner to be invertible. Finally, for this example, it is not necessary to impose any constraint since the projection on $Y_m$ leads to the best preconditioner and this preconditioner appears to be invertible for any $\xi \in \Xi$.

\subsubsection{Using the Frobenius semi-norm}

\begin{figure}[h!]
   \centering
   \subfigure[Interpolation functions.]{\centering
     \includegraphics[width=.45\textwidth]{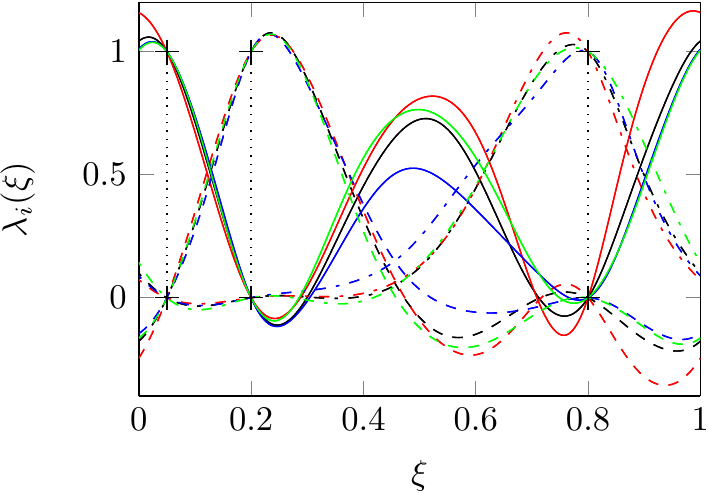}
     \label{fig:add_rot_Hadamard_MC_a}
     }~~
   \subfigure[Condition number of $P_3(\xi)A(\xi)$.]{\centering
     \includegraphics[width=.45\textwidth]{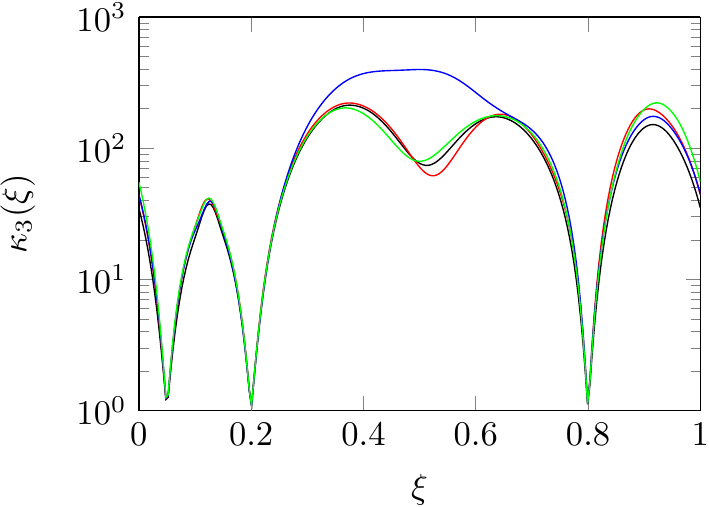}
     \label{fig:add_rot_Hadamard_MC_b}
     }
   \caption{Comparison between the Frobenius norm projection on $Y_3$ (black lines) and the Frobenius semi-norm projection on $Y_3$, using for $V$ either a sample of the rescaled Rademacher matrix (blue lines), the rescaled partial Hadamard matrix (red lines) or  a sample of the P-SRHT matrix (green lines) with $K=8$.}
   \label{fig:add_rot_Hadamard_MC}
\end{figure}

We analyze now the interpolation method defined by the Frobenius semi-norm projection on $Y_m$ \eqref{eq:SemiFrobProjection} for the different definitions of $V\in \Rbb^{n\times K}$ proposed in sections \ref{sec:Stat_trace} and \ref{sec:Hadamard_trace}. According to Table \ref{tab:error_lambda}, the error on the interpolation functions decreases slowly with $K$  (roughly as $\mathcal{O}( K^{-1/2} )$), and the use of the P-SRHT matrix leads to a slightly lower error. The interpolation functions are plotted on Figure \ref{fig:add_rot_Hadamard_MC_a} in the case where $K=8$. Even if we have an error of $36\%$ to $101\%$ on the interpolation functions, the condition number given on Figure \ref{fig:add_rot_Hadamard_MC_b} remains close to the one computed with the Frobenius norm. Also, an important remark is that with $K=8$ the computational effort for computing $M^V(\xi^*_k)$ and $S^V(\widetilde \xi^*_k)$ is negligible compared to the one for $M(\xi^*_k)$ and $S(\widetilde \xi^*_k)$.

\begin{table}[h!]
  \footnotesize
  \centering
  \begin{tabular}{|l|c|c|c|c|c|c|c|c|}
  \hline ~~~~~~~~~~~~~~~~~~$K$ & $8    $  & $16   $  & $ 32  $  & $64   $  & $128  $  & $256  $  & $512  $  \\ \hline
  Rescaled partial Hadamard & $0.4131$ & $0.3918$ & $0.3221$ & $0.1010$ & $0.0573$ & $0.0181$ & $0.0255$ \\ \hline
  Rescaled Rademacher (1)   & $0.5518$ & $0.0973$ & $0.2031$ & $0.1046$ & $0.1224$ & $0.1111$ & $0.0596$ \\
  Rescaled Rademacher (2)   & $1.0120$ & $0.6480$ & $0.1683$ & $0.1239$ & $0.0597$ & $0.0989$ & $0.0514$ \\
  Rescaled Rademacher (3)   & $0.7193$ & $0.2014$ & $0.1241$ & $0.1051$ & $0.1235$ & $0.1369$ & $0.0519$ \\ \hline
  P-SRHT (1)                & $0.4343$ & $0.2081$ & $0.2297$ & $0.0741$ & $0.0723$ & $0.0669$ & $0.0114$ \\
  P-SRHT (2)                & $0.3624$ & $0.2753$ & $0.0931$ & $0.1285$ & $0.0622$ & $0.0619$ & $0.0249$ \\
  P-SRHT (3)                & $0.8133$ & $0.4227$ & $0.1138$ & $0.0741$ & $0.0824$ & $0.0469$ & $0.0197$ \\\hline
  \end{tabular}
  \caption{Relative error $\sup_\xi \| \lambda(\xi) - \lambda^V(\xi) \|_{\mathbb{R}^3} / \sup_\xi \| \lambda(\xi) \|_{\mathbb{R}^3}$: $\lambda^V(\xi)$ (resp. $\lambda(\xi)$) are the interpolation functions associated to the Frobenius semi-norm projection (resp. the Frobenius norm projection) on $Y_m$, with $V$ either the rescaled partial Hadamard matrix, the random rescaled Rademacher matrix or the P-SRHT matrix (3 different samples for random matrices).}
  \label{tab:error_lambda}
\end{table}

\subsubsection{Greedy selection of the interpolation points}

\begin{figure}[h!]
   \centering
   \subfigure{\centering
     \includegraphics[width=.31\textwidth]{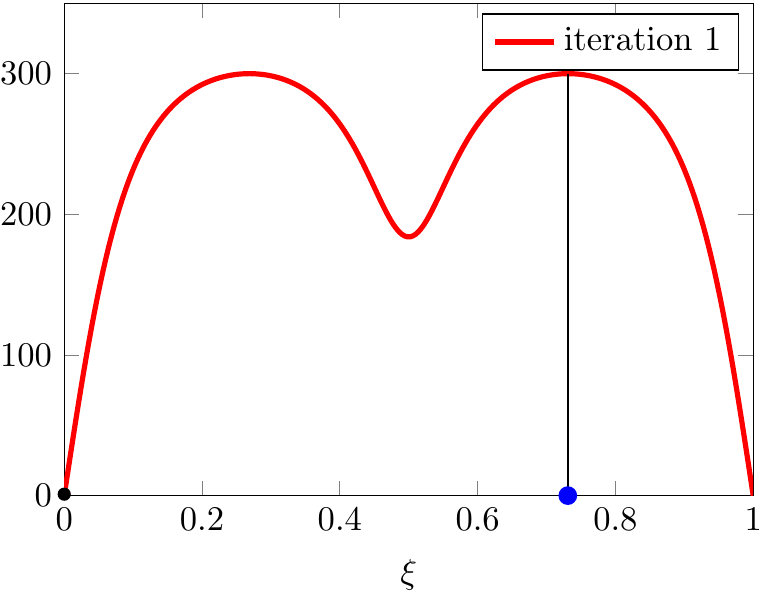}
     }
   \subfigure{\centering
     \includegraphics[width=.31\textwidth]{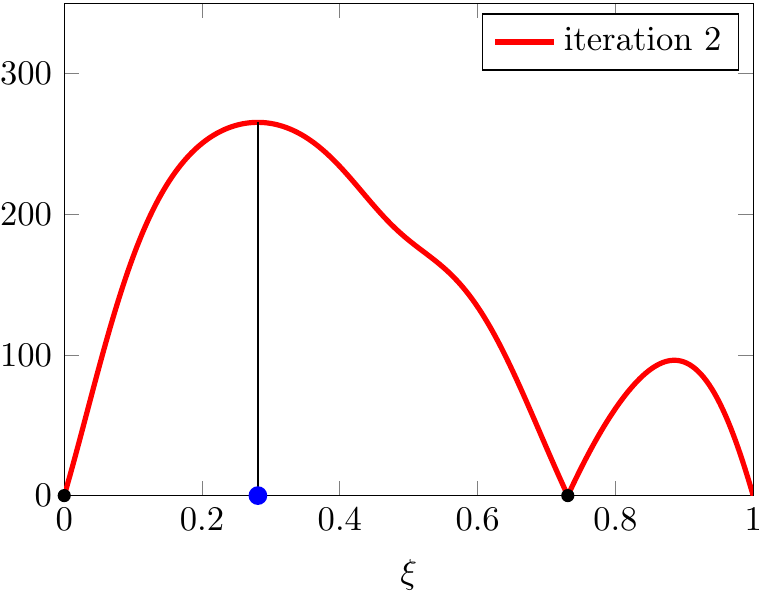}
     }
   \subfigure{\centering
     \includegraphics[width=.31\textwidth]{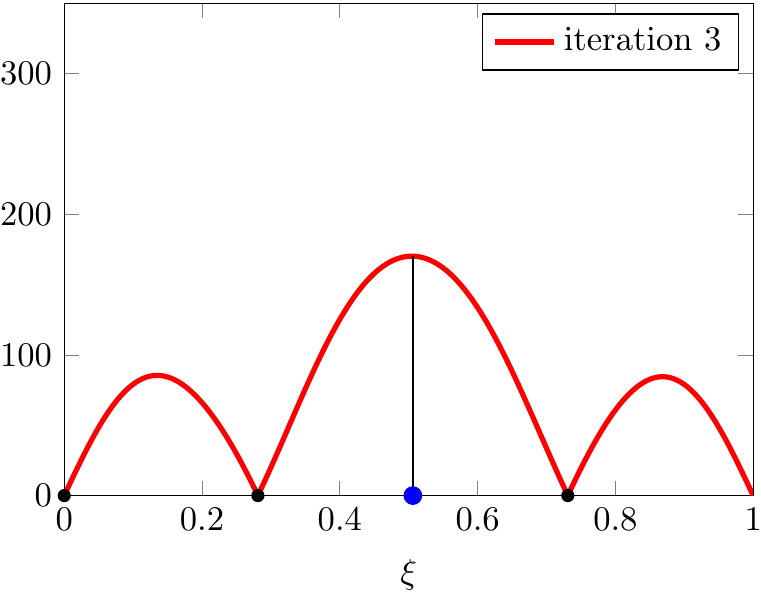}
     }
   \subfigure{\centering
     \includegraphics[width=.31\textwidth]{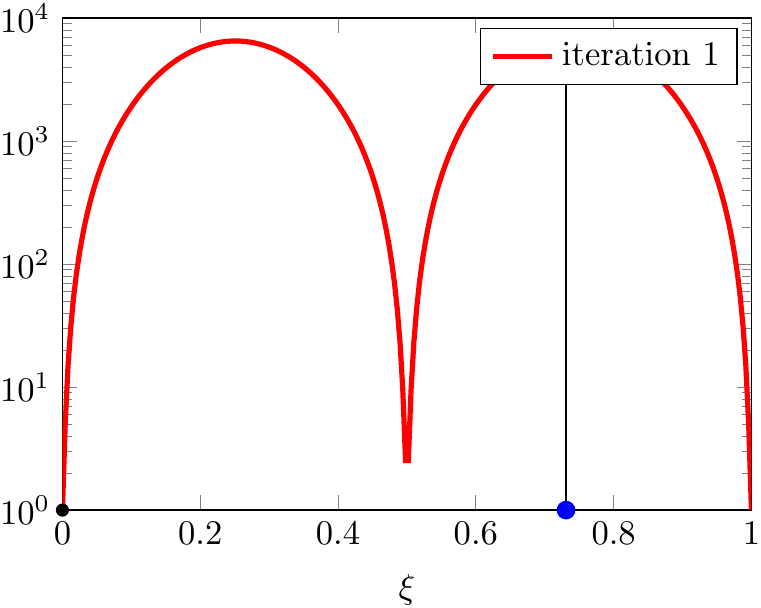}
     \label{fig:add_rot_greedy_a}
     }~~
   \subfigure{\centering
     \includegraphics[width=.31\textwidth]{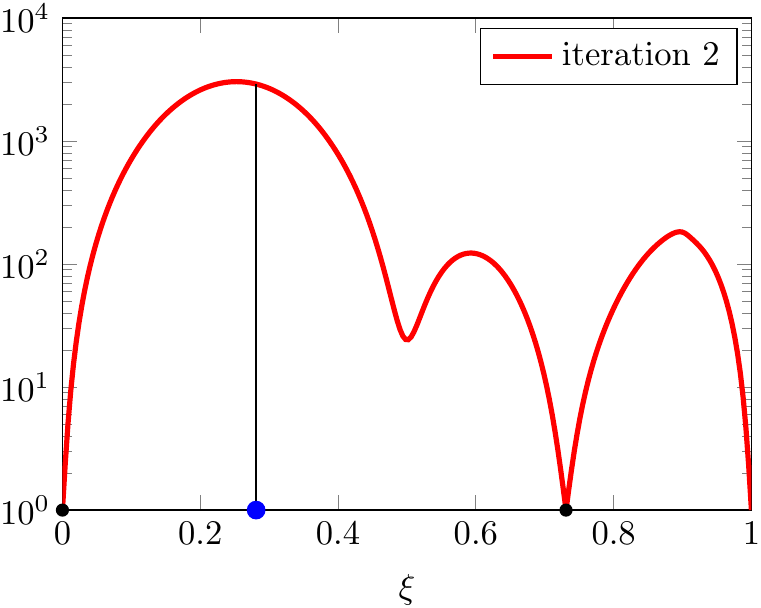}
     \label{fig:add_rot_greedy_b}
     }
   \subfigure{\centering
     \includegraphics[width=.31\textwidth]{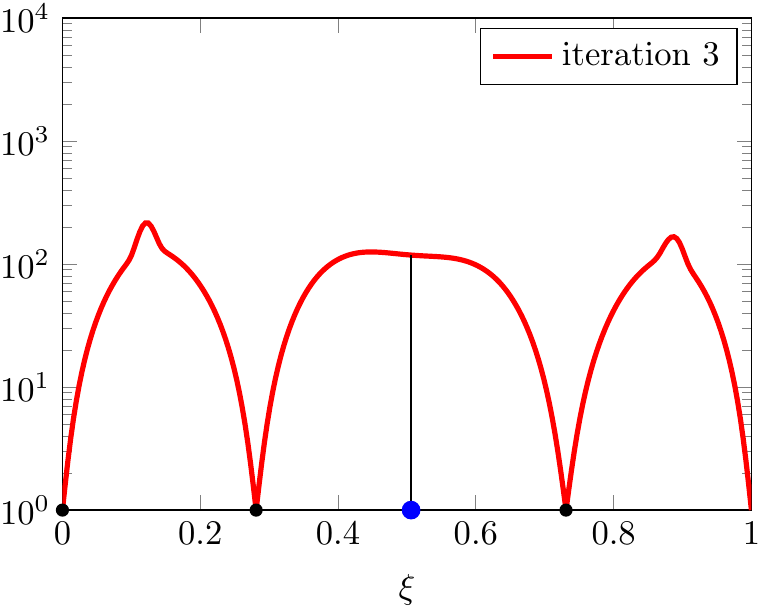}
     \label{fig:add_rot_greedy_c}
     }

   \caption{Greedy selection of the interpolation points: the first row is the residual $\| (I-P_k(\xi)A(\xi)) V\|_F$ (the blue points correspond to the maximum of the residual) with $V$ a realization of the P-SRHT matrix with $K=128$ columns, and the second row is the condition number of $P_m(\xi)A(\xi)$.}
   \label{fig:add_rot_greedy}
\end{figure}

We now consider the greedy selection of the interpolation points presented in Section \ref{sec:Greedy_Precond}. We start with an initial point $\xi_1=0$ and the next points are defined by \eqref{eq:greedy_precond_frob}, where matrix $V$ is a realization of the P-SRHT matrix with $K=128$ columns. $P_m(\xi)$ is the projection on $Y_m$ using the Frobenius semi-norm defined by \eqref{eq:SemiFrobProjection}. The first 3 steps of the algorithm are illustrated on Figure \ref{fig:add_rot_greedy}. We observe that at each iteration, the new interpolation point $\xi_{m+1}$  is close to the point where the condition number of $P_m(\xi)A(\xi)$ is maximal. Table \ref{tab:add_rot_greedy} presents the maximal value over $\xi\in\Xi$ of the residual, and of the condition number of $P_m(\xi)A(\xi)$. Both quantities are rapidly decreasing with  $m$. This shows that this algorithm, initially designed to minimize $\| (I-P_m(\xi)A(\xi) )V \|_F$, seems to be also efficient for the construction of preconditioners, in the sense that the condition number decreases rapidly.

\begin{table}[h!]
  \footnotesize
  \centering
  \begin{tabular}{|c|c|c|c|c|c|c|c|}
  \hline iteration $m$                      & 0      &  1     &  2     &  5      &  10     &  20   &  30    \\ \hline
  $\sup_\xi \kappa(P_m(\xi)A(\xi))$         & 10001  &  6501  &  3037  &  165,7  &  51,6   &  16,7 &  7,3   \\
  $\sup_\xi \| (I-P_m(\xi)A(\xi)) V\|_F$    & -      &  300   &  265   &  80,5   &  35,4   &  10,0 &  7,6   \\ \hline
  \end{tabular}
  \caption{Convergence of the greedy algorithm: supremum over $\xi\in\Xi$ of the condition number (first row) and of the Frobenius semi-norm residual (second row).}
  \label{tab:add_rot_greedy}
\end{table}

\subsection{Multi-parameter-dependent equation}

\begin{figure}[h!]
   \centering
   \subfigure{\centering
     \includegraphics[scale=0.65]{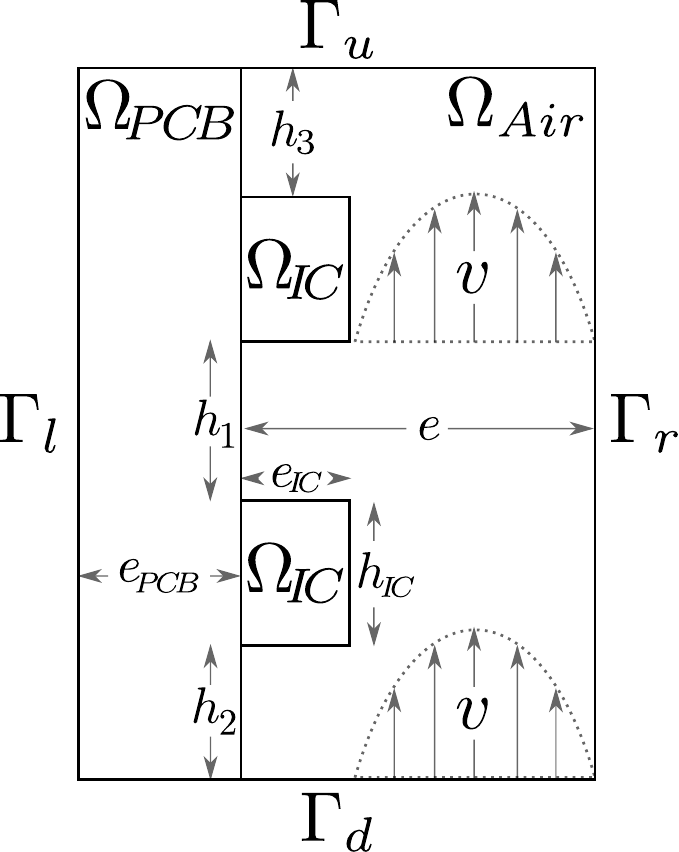}
     \label{fig:OPUS_geometry}
     }
   \subfigure{\centering
     \includegraphics[scale=0.8]{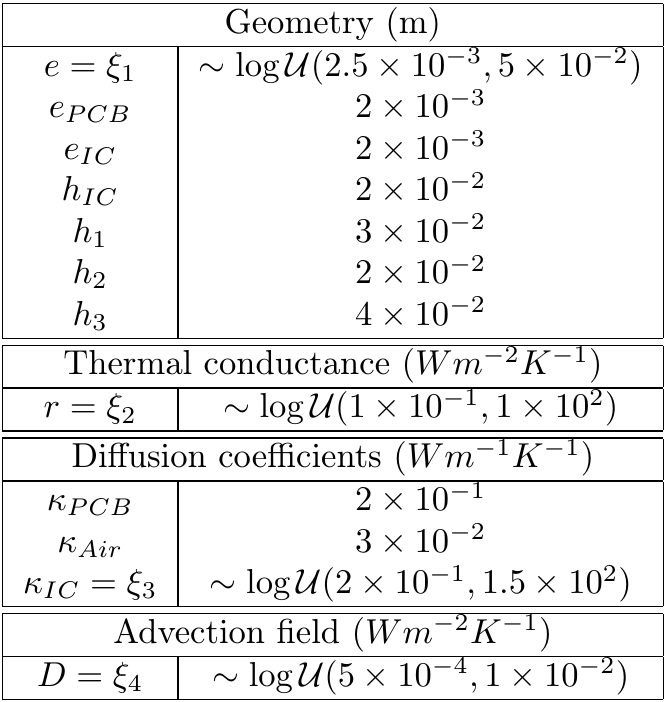}
     \label{fig:OPUS_data}
     }
   \caption{Geometry and parameters of the benchmark OPUS.}
   \label{fig:Opus_settings}
\end{figure}

We introduce a benchmark proposed within the OPUS project (see http://www.opus-project.fr). Two electronic components $\Omega_{IC}$ (see Figure \ref{fig:Opus_settings}) submitted to a cooling air flow in the domain $\Omega_{Air}$ are fixed on a printed circuit board $\Omega_{PCB}$. The temperature field defined over $\Omega=\Omega_{IC}\cup \Omega_{PBC}\cup \Omega_{Air}\subset \mathbb{R}^2$ satisfies the advection-diffusion equation:
\begin{align}\label{eq:OPUS}
 -\nabla\cdot( \kappa(\xi) \nabla u ) +D(\xi)v\cdot\nabla u &= f.
\end{align}

The diffusion coefficient $\kappa(\xi)$ is equal to $\kappa_{{PCB}}$ on $\Omega_{PCB}$, $\kappa_{Air}$ on  $\Omega_{Air}$ and $\kappa_{IC}$ on $\Omega_{IC}$. The right hand side $f$ is equal to $Q=10^{6}$ on $\Omega_{IC}$ and $0$ elsewhere. The boundary conditions are $u=0$ on $\Gamma_d$, $e_2\cdot\nabla u=0$ on $\Gamma_u$ ($e_1,e_2$ are the canonical vectors of $\mathbb{R}^2$), and $u_{|\Gamma_l} = u_{|\Gamma_r}$ (periodic boundary condition).
At the interface $\Gamma_C=\partial \Omega_{IC}\cap \partial\Omega_{PCB}$ there is a thermal contact conductance, meaning that the temperature field $u$ admits a jump over $\Gamma_C$ which satisfies
\begin{equation*}
 \kappa_{IC}(e_1\cdot\nabla u_{|\Omega_{IC}}) = \kappa_{PCB}(e_1\cdot\nabla u_{|\Omega_{PCB}})=r(u_{|\Omega_{IC}}-u_{|\Omega_{PCB}})~~~\mbox{ on } \Gamma_C.
\end{equation*}

The advection field $v$ is given by $v(x,y)=e_2 g(x) $, where $g(x)=0$ if $x\leq e_{PCB}+e_{IC}$ and 
\begin{equation*}
 g(x) = \frac{3}{2(e-e_{IC})} \left( 1-\left( \frac{2x- (e+e_{IC}+2e_{PCB}) }{e-e_{IC}}\right)^2 \right)
\end{equation*}
otherwise. 
We have 4 parameters: the width $e:=\xi_1$ of the domain $\Omega_{Air}$, the thermal conductance parameter $r:=\xi_2$, the diffusion coefficient  $\kappa_{IC}:=\xi_3$ of the components and the amplitude of the advection field $D:=\xi_4$. Since the domain $\Omega=\Omega(e)$ depends on the parameter $\xi_1\in[e_{min},e_{max}]$, we introduce a geometric transformation $(x,y)=\phi_{\xi_1}(x_0,y_0)$ that maps a reference domain $\Omega_0=\Omega(e_{max})$ to $\Omega(\xi_1)$:
\begin{equation*}
 \phi_{\xi_1}(x_0,y_0) = \begin{pmatrix}
                          \left\{ \begin{matrix} x_0  &\mbox{ if }x_0\leq e_0 \\ e_0 + (x_0-e_0)\frac{\xi_1 - e_{IC}}{e_{max} - e_{IC}} &\mbox{ otherwise.}\end{matrix}\right\}
			  \\y_0
                         \end{pmatrix} ,
\end{equation*}
with $e_0 = e_{PCB}+e_{IC}$. This method is described in \cite{Rozza2008}: since the geometric transformation $\phi_{\xi_1}$ satisfies the so-called \textit{Affine Geometry Precondition}, the operator of equation \eqref{eq:OPUS} formulated on the reference domain admits an affine representation.

For the spatial discretization we use a finite element approximation with $n=2.8\cdot 10^{4}$ degrees of freedom (piecewise linear approximation). We rely on a Galerkin method with SUPG stabilization (see \cite{Brooks1982}). $\Xi$ is a set of $10^4$ independent samples drawn according the loguniform probability laws of the parameters given on Figure \ref{fig:Opus_settings}.

\subsubsection{Preconditioner for the projection on a given reduced space} We consider here a POD basis $X_r$ of dimension $r=50$ computed with $100$ snapshots of the solution (a basis of $X_r$ is obtained by the first $50$ dominant singular vectors of a matrix of 100 random snapshots of $u(\xi)$). Then  we compute the Petrov-Galerkin projection as presented in Section \ref{sec:MR_optimal_proj}. The efficiency of the preconditioner can be measured with the quantity $\delta_{r,m}(\xi)$: the associated quasi-optimality constant $(1-\delta_{r,m}(\xi)^2)^{-1/2}$ should be as close to one as possible (see equation \eqref{eq:control_PG_2}). We introduce the quantile $q_p$ of probability $p$ associated to the quasi-optimality constant $(1-\delta_{r,m}(\xi)^2)^{-1/2}$ defined as the smallest value $q_p\geq 1$ satisfying
\begin{equation*}
 \mathbb{P} \big( \{\xi\in\Xi :  (1-\delta_{r,m}(\xi)^2 )^{-1/2}\leq q_p \} \big) \geq p,
\end{equation*}
where $\mathbb{P}(A)= \#A/\#\Xi$ for $A\subset \Xi$.
Table \ref{tab:opus_delta} shows the evolution of the quantile with respect to the number of interpolation points for the preconditioner. Here the goal is to compare the different strategies for the selection of the interpolation points:
\begin{itemize}
 \item[(a)] the greedy selection \eqref{eq:greedy_precond_delta} based on the quantity $\delta_{r,m}(\xi)$, 
 \item[(b)] the greedy selection \eqref{eq:greedy_precond_frob} based on the Frobenius semi-norm residual, with $V$ a P-SRHT matrix 
with $K=256$ columns, and 
 \item[(c)] a random Latin Hypercube sample (LHS).
\end{itemize}
 The projection on $Y_m$ (or $Y_m^+$) is then defined with the Frobenius semi-norm using for $V$ a P-SRHT matrix with $K=330$ columns.

When considering a small number of interpolation points $m\leq3$, the projection on $Y_m^+$ provides lower quantiles for the quasi-optimality constant compared to the projection on $Y_m$. The positivity constraint is useful for small $m$. But for high values of $m$ (see $m=15$) the positivity constraint is no longer necessary and the projection on $Y_m$ provides lower quantiles.

Concerning the choice of the interpolation points, the strategy (a) shows the faster decay of the quantiles $q_p$, especially for $p=50\%$. The strategy (b) shows also good results, but the quantile $q_p$ for $p=100\%$ are still high compared to (a). These results show the benefits of the greedy selection based on the quasi-optimality constant. Finally the strategy (c) shows bad results (high values of the quantiles), especially for small $m$.

\begin{table}[h!]
  \footnotesize
  \centering
  \begin{tabular}{c|c|c|c|c|c|c|c|c|c|}
    \cline{2-10}
      & \multicolumn{9}{c|}{Projection on $Y_m$} \\ 
      & \multicolumn{6}{|c|}{Greedy selection based on} & \multicolumn{3}{c|}{(c) Latin Hypercube} \\ 
      & \multicolumn{3}{|c|}{(a) $\delta_{r,m}(\xi)$} &  \multicolumn{3}{ c| }{(b) Frob. residual} & \multicolumn{3}{c|}{sampling} \\ 
      &$50\%$&$90\%$&$100\%$&$50\%$&$90\%$&$100\%$&$50\%$&$90\%$&$100\%$ \\ \cline{1-10}
      \multicolumn{1}{|l|}{$m=0$}  & $ 21.3 $ & $ 64.1 $ & $ 94.1 $  & $ 21.3 $ & $ 64.1 $ & $ 94.1 $  & $ 21.3 $ & $ 64.1 $ & $ 94.1 $ \\
      \multicolumn{1}{|l|}{$m=1$}  & $ 18.3 $ & $ 74.1 $ & $ 286.7 $ & $ 10.2 $ & $ 36.1 $ & $ 161.6 $ & $ 18.3 $ & $ 104.1 $ & $ 231.8 $ \\
      \multicolumn{1}{|l|}{$m=2$}  & $ 11.9 $ & $ 22.6 $ & $ 42.1 $  & $ 9.8 $  & $ 53.3 $ & $ 374.0 $ & $ 11.5 $ & $ 113.0 $ & $ 533.9 $ \\
      \multicolumn{1}{|l|}{$m=3$}  & $ 11.1 $ & $ 49.2 $ & $ 200.4 $ & $ 7.8 $  & $ 31.2 $ & $ 60.2 $  & $ 18.3 $ & $ 138.7 $ & $ 738.5 $ \\
      \multicolumn{1}{|l|}{$m=5$}  & $ 5.2 $  & $ 10.8 $ & $ 18.4 $  & $ 6.8 $  & $ 18.6 $ & $ 24.5 $  & $ 8.7 $ & $ 121.1 $ & $ 651.4 $ \\
      \multicolumn{1}{|l|}{$m=10$} & $ 3.1 $  & $ 9.0 $  & $ 13.2 $  & $ 5.3 $  & $ 22.3 $ & $ 62.1 $  & $ 4.0 $ & $ 21.6 $ & $ 345.7 $ \\
      \multicolumn{1}{|l|}{$m=15$} & $ 2.2 $  & $ 6.3 $  & $ 10.4 $  & $ 3.5 $  & $ 6.5 $  & $ 11.5 $  & $ 2.7 $ & $ 7.8 $ & $ 48.6 $ 
      \\ \hline \multicolumn{10}{c}{ } \\
     \cline{2-10}
      & \multicolumn{9}{c|}{Projection on $Y_m^+$} \\ 
      & \multicolumn{6}{|c|}{Greedy selection based on} & \multicolumn{3}{c|}{(c) Latin Hypercube} \\ 
      & \multicolumn{3}{|c|}{(a) $\delta_{r,m}(\xi)$} &  \multicolumn{3}{ c| }{(b) Frob. residual} & \multicolumn{3}{c|}{sampling} \\ 
      &$50\%$&$90\%$&$100\%$&$50\%$&$90\%$&$100\%$&$50\%$&$90\%$&$100\%$ \\ \cline{1-10}
      \multicolumn{1}{|l|}{$m=0$}  & $ 21.3 $ & $ 64.1 $ & $ 94.1 $  & $ 21.3 $ & $ 64.1 $ & $ 94.1 $  & $ 21.3 $ & $ 64.1 $  & $ 94.1 $ \\
      \multicolumn{1}{|l|}{$m=1$}  & $ 18.3 $ & $ 74.1 $ & $ 286.7 $ & $ 10.2 $ & $ 36.1 $ & $ 161.6 $ & $ 18.3 $ & $ 104.1 $ & $ 231.8 $\\
      \multicolumn{1}{|l|}{$m=2$}  & $ 11.9 $ & $ 22.6 $ & $ 42.1 $  & $ 8.9 $  & $ 35.5 $ & $ 78.6 $  & $ 10.4 $ & $ 41.5 $  & $ 112.5 $\\
      \multicolumn{1}{|l|}{$m=3$}  & $ 9.7 $  & $ 24.4 $ & $ 48.0 $  & $ 7.9 $  & $ 27.7 $ & $ 57.9 $  & $ 12.1 $ & $ 48.8 $  & $ 114.1 $\\
      \multicolumn{1}{|l|}{$m=5$}  & $ 6.4 $  & $ 15.0 $ & $ 25.5 $  & $ 6.9 $  & $ 26.8 $ & $ 65.1 $  & $ 5.7 $  & $ 11.6 $  & $ 17.5 $ \\
      \multicolumn{1}{|l|}{$m=10$} & $ 4.6 $  & $ 9.5 $  & $ 16.8 $  & $ 7.3 $  & $ 18.9 $ & $ 38.0 $  & $ 4.3 $  & $ 10.0 $  & $ 18.5 $ \\
      \multicolumn{1}{|l|}{$m=15$} & $ 4.3 $  & $ 7.1 $  & $ 11.2 $  & $ 6.4 $  & $ 10.1 $ & $ 18.0 $  & $ 4.2 $  & $ 9.0 $   & $ 19.3 $ \\ \hline
  \end{tabular}
  \caption{Quantiles $q_p$ of the quasi-optimality constant associated to the Petrov-Galerkin projection on the POD subspace $X_r$ for $p=50\%,$ $90\%$ and $100\%$. The row $m=0$ corresponds to $P_0(\xi)=R_X^{-1}$, that is the standard Galerkin projection.}
  \label{tab:opus_delta}
\end{table}

\subsubsection{Preconditioner for Reduced Basis method} 
\begin{figure}[h!]
   \centering
   \subfigure[]{\centering
     \includegraphics[scale=0.25]{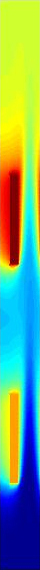}
     \label{fig:OPUS_sampling_a}
     }~~~~~~
     \subfigure[]{\centering
     \includegraphics[scale=0.25]{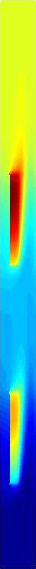}
     \label{fig:OPUS_sampling_b}
     }~~~~~~
     \subfigure[]{\centering
     \includegraphics[scale=0.25]{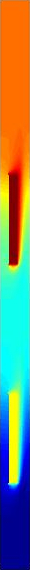}
     \label{fig:OPUS_sampling_c}
     }~~~~~~
     \subfigure[]{\centering
     \includegraphics[scale=0.25]{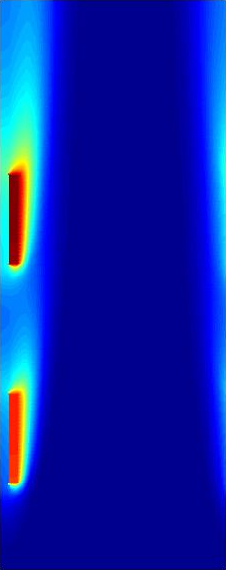}
     \label{fig:OPUS_sampling_d}
     }~~~~~~
     \subfigure[]{\centering
     \includegraphics[scale=0.25]{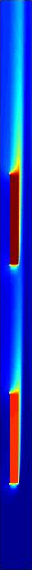}
     \label{fig:OPUS_sampling_e}
     }~~~~~~
     \subfigure[]{\centering
     \includegraphics[scale=0.25]{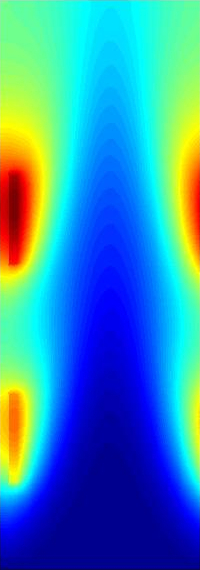}
     \label{fig:OPUS_sampling_f}
     }\\
     {\footnotesize \begin{tabular}{r|l|l|l|l|l|l}
                    & (a) $r=1$          & (b) $r=2$          & (c) $r=3$          & (d)  $r=4$         & (e) $r=5$ & (f) $r=6$ \\ \hline
      $e$           & $7.1\cdot 10^{-3}$ & $6.1\cdot 10^{-3}$ & $5.0\cdot 10^{-3}$ & $5.0\cdot 10^{-2}$ & $5.1\cdot 10^{-3}$ & $4.4\cdot 10^{-2}$\\ 
      $r$           & $9.9\cdot 10^{1 }$ & $2.3\cdot 10^{0 }$ & $4.1\cdot 10^{-1}$ & $2.8\cdot 10^{0 }$ & $8.6\cdot 10^{-1}$ & $4.8\cdot 10^{1 }$\\ 
      $\kappa_{IC}$ & $1.1\cdot 10^{2 }$ & $2.1\cdot 10^{-1}$ & $3.2\cdot 10^{1 }$ & $6.0\cdot 10^{0 }$ & $1.1\cdot 10^{2 }$ & $3.2\cdot 10^{-1}$\\ 
      $D$           & $1.7\cdot 10^{-3}$ & $7.4\cdot 10^{-4}$ & $6.2\cdot 10^{-4}$ & $9.9\cdot 10^{-3}$ & $8.6\cdot 10^{-3}$ & $6.0\cdot 10^{-4}$
     \end{tabular}
     }
   \caption{First six interpolation points of the ideal reduced basis method and corresponding reduced basis functions.}
   \label{fig:OPUS_sampling}
\end{figure}

\begin{figure}[h!]
   \centering
   \subfigure[Projection on $Y_m$]{\centering
     \includegraphics[width=.45\textwidth]{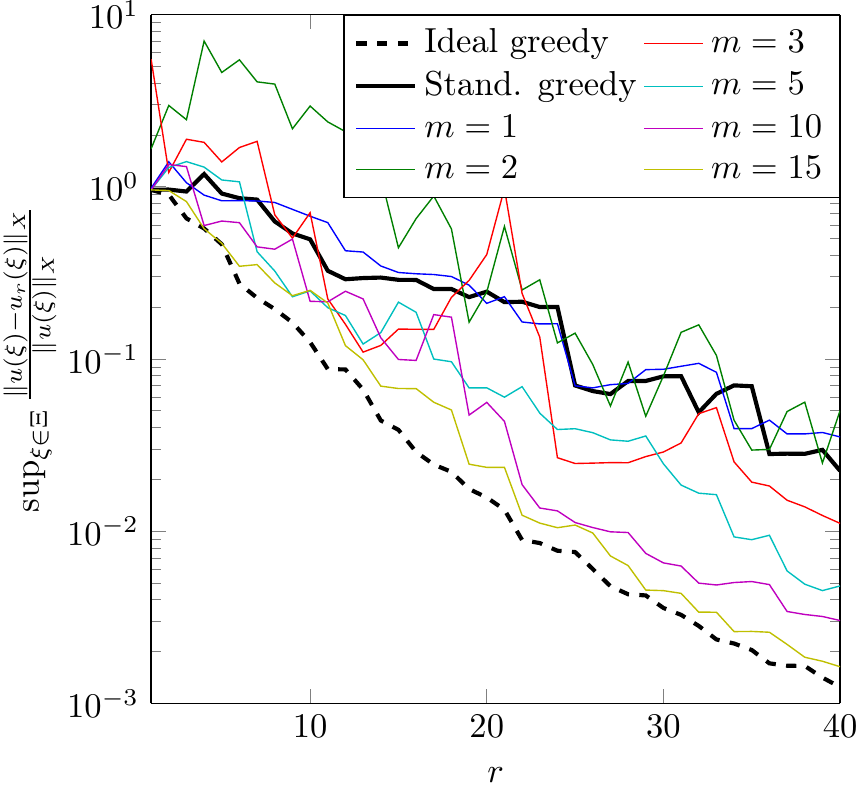}
     \label{fig:OPUS_RB_Greedy_SC}
     }~
   \subfigure[Projection on $Y_m^+$]{\centering
     \includegraphics[width=.45\textwidth]{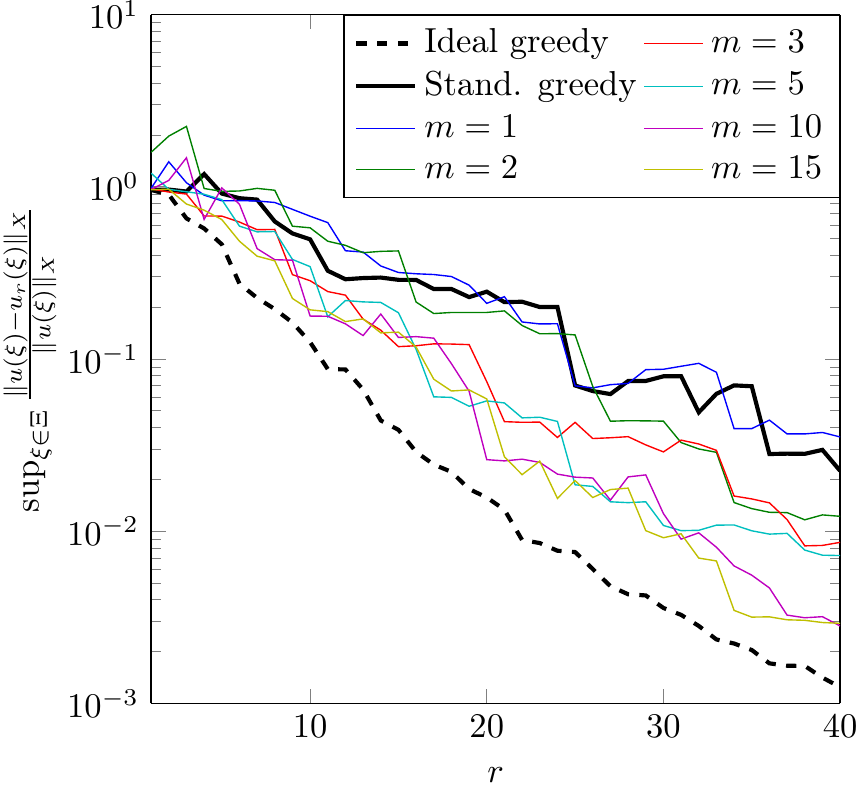}
     \label{fig:OPUS_RB_Greedy_AC}
     } \\
   \subfigure[Projection on $Y_m$]{\centering
     \includegraphics[width=.45\textwidth]{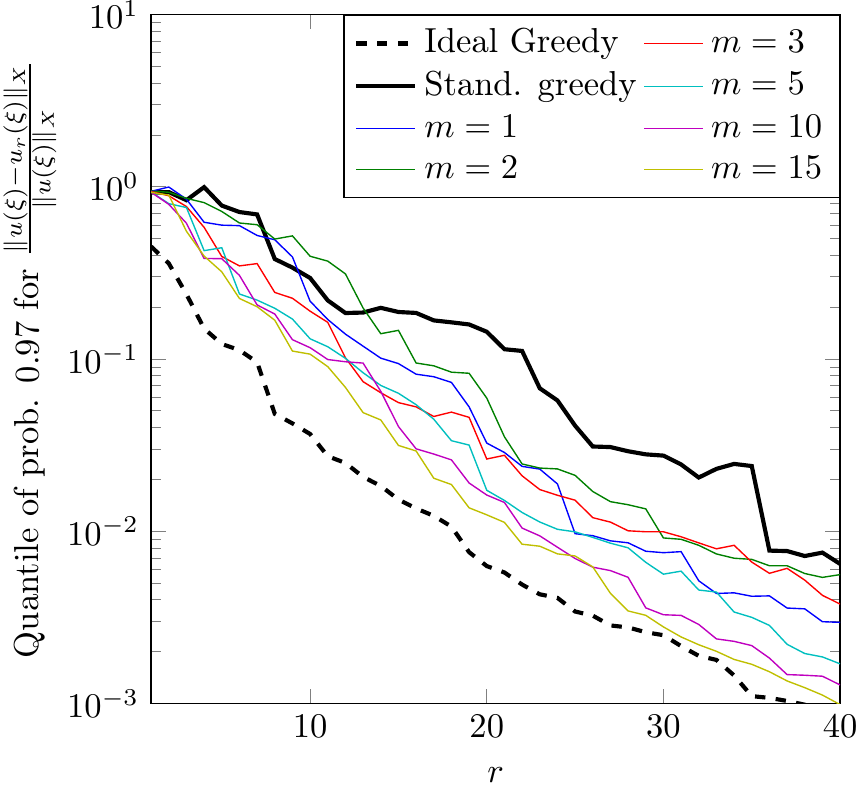}
     \label{fig:OPUS_RB_Greedy_SC_quantile}
     }~
   \subfigure[Projection on $Y_m^+$]{\centering
     \includegraphics[width=.45\textwidth]{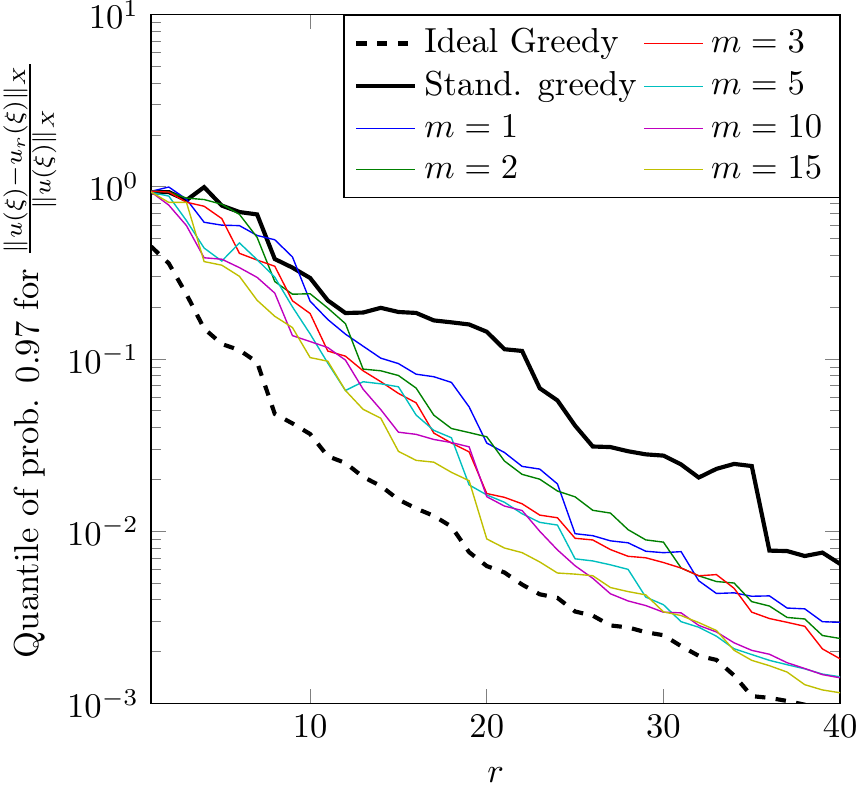}
     \label{fig:OPUS_RB_Greedy_AC_quantile}
     }
   \caption{Convergence of the preconditioned reduced basis method using the greedy selection of interpolation points for the preconditioner. Supremum over $\Xi$ (top) and quantile of probability $97\%$ (bottom) of the relative error $\| u(\xi)-u_r(\xi) \|_X/\| u(\xi)\|_X$ with respect to $r$. Comparison of preconditioned reduced basis algorithms with ideal and standard greedy algorithms.}
   \label{fig:OPUS_reduced_basis_GREEDY}
\end{figure}

\begin{figure}[h!]
  \centering
  \includegraphics[width=.45\textwidth]{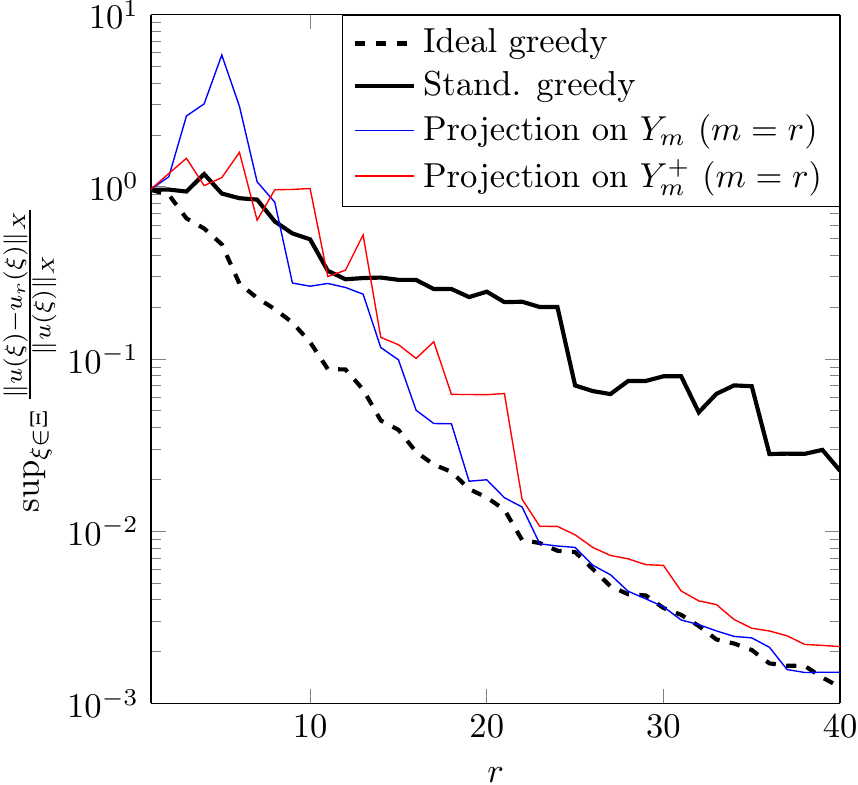}
  \includegraphics[width=.45\textwidth]{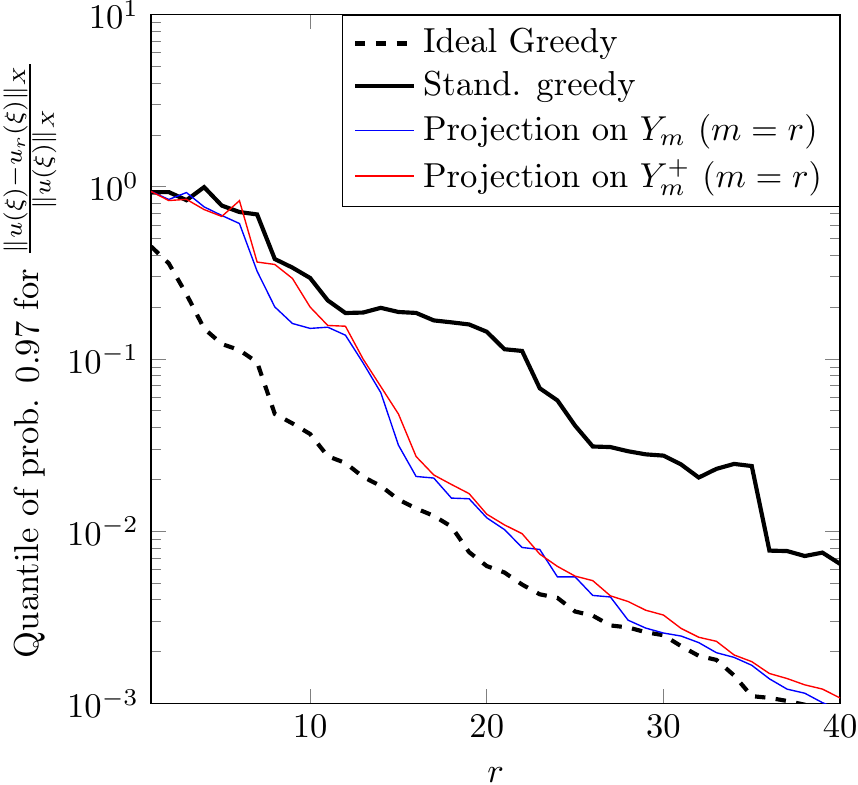}
  \caption{Preconditioned Reduced basis methods with re-use of operators. Supremum over $\Xi$ (left) and quantile of probability $97\%$ (right) of the relative error $\| u(\xi)-u_r(\xi) \|_X/\| u(\xi)\|_X$ with respect to $r$. Comparison of preconditioned reduced basis algorithms with ideal and standard  greedy algorithms.}
  \label{fig:OPUS_reduced_basis_recycling}
\end{figure}

We now consider the preconditioned Reduced Basis method for the construction of the approximation space $X_r$, as presented in Section \ref{sec:MR_greedy_approx}. 
Figures \ref{fig:OPUS_reduced_basis_GREEDY} and \ref{fig:OPUS_reduced_basis_recycling} show the convergence of the error with respect to the rank $r$ of $u_r(\xi)$ for different constructions of the preconditioner $P_m(\xi)$. Two measures of the error are given: $\sup_{\xi\in\Xi} \| u(\xi)-u_r(\xi) \|_X/\| u(\xi) \|_X$, and the quantile of probability $0.97$ for $\| u(\xi)-u_r(\xi) \|_X/\| u(\xi) \|_X$. The curve ``Ideal greedy'' corresponds to the algorithm defined by \eqref{eq:ideal_greedy} which provides a reference for the ideally conditioned algorithm, \textit{i.e.} with $\kappa_m(\xi)=1$. Figure \ref{fig:OPUS_sampling} shows the corresponding first interpolation points for the solution.

The greedy selection of the interpolation points based on \eqref{eq:greedy_precond_frob} (see Figure \ref{fig:OPUS_reduced_basis_GREEDY}) allows to almost recover the convergence curve of the ideal greedy algorithm when using the projection on $Y_m$ with $m=15$. For the strategy of re-using the operators factorizations, the approximation is rather bad for $r=m\leq 10$ meaning that the space $Y_r$ (or $Y_r^+$) is not really adapted for the construction of a good preconditioner over the whole parametric domain. However, for higher values of $r$, the preconditioner is getting better and better. For $r\geq 20$, we almost reach the convergence of the ideal greedy algorithm. We conclude that this strategy of re-using the operator factorization, which has a computational cost comparable to the standard non preconditioned Reduced Basis greedy algorithm, allows obtaining asymptotically the performance of the ideal greedy algorithm. Note that the positivity constraint yields a better preconditioner for small values of $r$ but is no longer necessary for large $r$.

Let us finally consider the effectivity index 
\begin{equation*}
 \eta_r(\xi)  = {\| P_r(\xi) (A(\xi)u_r(\xi) - b(\xi) )\|_{X}}/{ \| u(\xi) - u_r(\xi) \|_X },
\end{equation*}
which evaluates the quality of the preconditioned residual norm for error estimation. We introduce the confidence interval $I_r(p)$ defined as the smallest interval which satisfies
\begin{equation*}
\mathbb{P}(\{\xi\in\Xi~:~ \eta_r(\xi)\in I_r(p) \}) \geq p.
\end{equation*}

\begin{figure}[h!]
  \centering
  \subfigure[Re-use, proj. on $Y_m$.]{\centering
     \includegraphics[width=.3\textwidth]{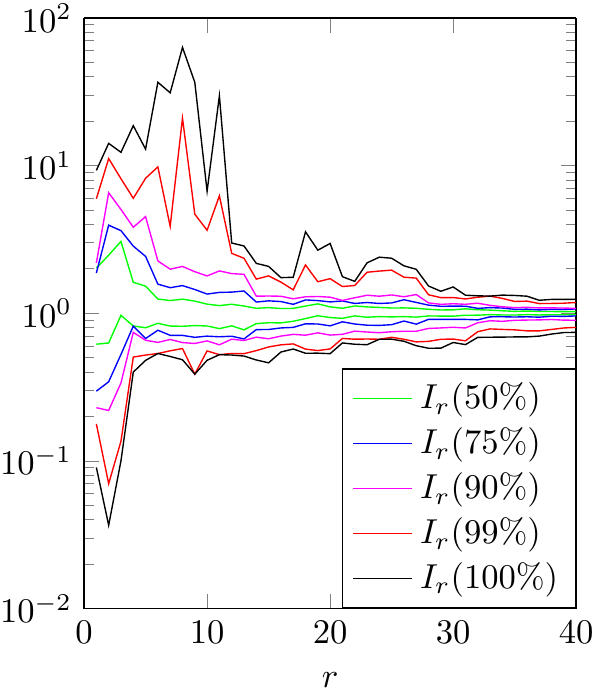}
     \label{fig:OPUS_R_SC_effectivity}
     }
  \subfigure[Re-use, proj. on $Y_m^+$.]{\centering
     \includegraphics[width=.3\textwidth]{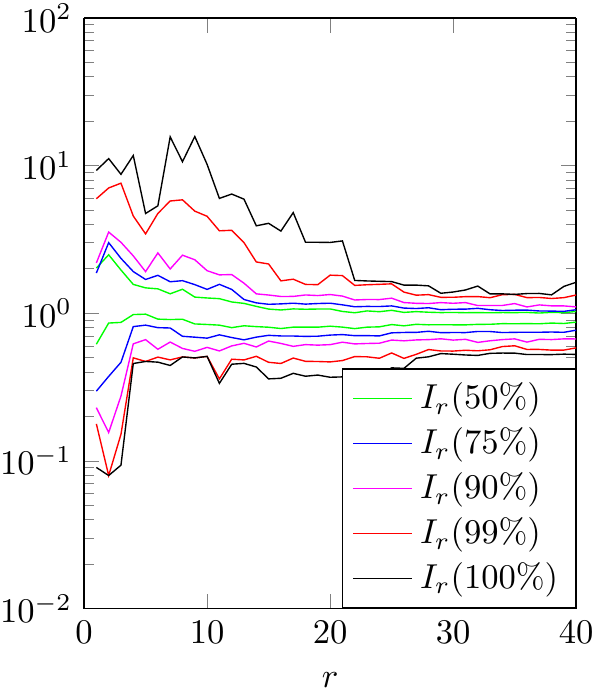}
     \label{fig:OPUS_R_AC_effectivity}
     }
  \subfigure[Dual residual norm.]{\centering
     \includegraphics[width=.3\textwidth]{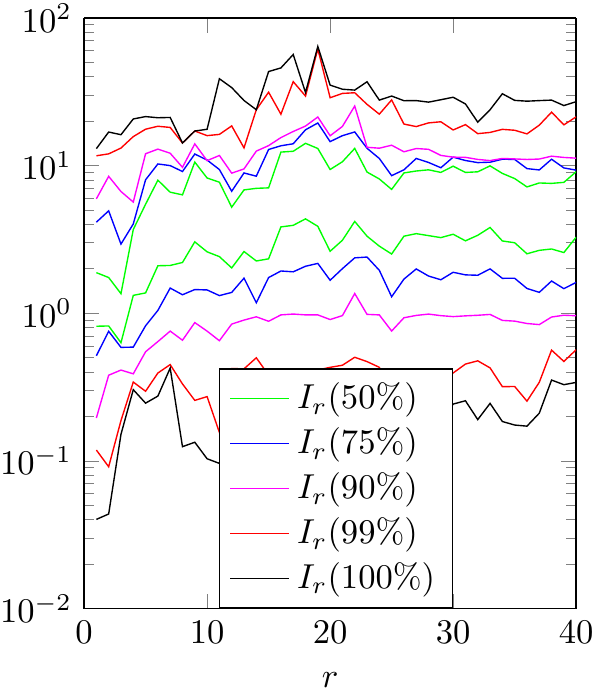}
     \label{fig:OPUS_CRB_effectivity}
     }
  \caption{Confidence intervals of the effectivity index during the iterations of the Reduced Basis greedy construction. Comparison between preconditioned algorithms with re-use of operators factorizations (a,b) and the non preconditioned greedy algorithm (c).}
  \label{fig:OPUS_effectivity}
\end{figure}

On Figure \ref{fig:OPUS_effectivity} we see that the confidence intervals are shrinking around $1$ when $r$ increases, meaning that the preconditioned residual norm becomes a better and better error estimator when $r$ increases. Again, the positivity constraint is needed for small values of $r$, but we obtain a better error estimation without imposing this constraint for $r\geq20$. On the contrary, the standard residual norm leads to effectivity indices that spread from $10^{-1}$ to $10^{1}$ with no improvement as $r$ increases, meaning that we can have a factor $10^2$ between the error estimator $\| A(\xi) u_r(\xi)-b(\xi) \|_{X'}$ and the true error $\| u_r(\xi) - u(\xi) \|_X$.

\section{Conclusion} We have proposed a method for the interpolation of the inverse of a parameter-dependent matrix. The interpolation is defined by the projection of the identity in the sense of the Frobenius norm. Approximations of the Frobenius norm have been introduced to make computationally feasible the projection in the case of large matrices. Then, we have proposed preconditioned versions of projection-based model reduction methods. The preconditioner can be used to define Petrov-Galerkin projections on a given approximation space with better quasi-optimality constants by introducing a parameter-dependent test space depending on the preconditioner. Also, the preconditioner can be used to improve residual-based error estimates that are used for assessing the quality of a given approximation, which is required in any adaptive approximation strategy. Different strategies have  been proposed for the selection of interpolation points depending on the objective: (i) the construction of an optimal approximation of the inverse operator for preconditioning iterative solvers or for improving error estimators based on preconditioned residuals, (ii) the improvement of the quality of Petrov-Galerkin projections of the solution of a parameter-dependent equation on a given reduced approximation space, or (iii) the re-use of operators factorizations when solving a parameter-dependent equation with a sample-based approach. The performance of the obtained parameter-dependent preconditioners has been illustrated in the context of projection-based model reduction techniques such as the Proper Orthogonal Decomposition and the Reduced Basis method.

The proposed preconditioner has been used to define Petrov-Galerkin projections with better stability constants. For the solution of PDEs, the Petrov-Galerkin projection has been defined at the discrete (algebraic) level for obtaining a better approximation (in a reduced space) of the finite element Galerkin approximation of the PDE. Therefore, for convection-dominated problems, the proposed approach does not avoid using stabilized finite element formulations. Similar observations can be found in \cite{Pacciarini2014}. However, a Petrov-Galerkin method could be defined at the continuous level with a preconditioner being the interpolation of inverse partial differential operators. In this continuous framework, the preconditioner would improve the stability constant for the finite element Galerkin projection and may avoid the use of stabilized finite element formulations. Such Petrov-Galerkin methods  have been proposed in \cite{Cohen:2012,Dahmen2011} for convection-dominated problems (as an alternative to standard stabilization methods), which can be interpreted as an implicit preconditioning method defined at the continuous level.

In the present paper, the parameter-dependent preconditioner is obtained by a projection onto the space generated by snapshots of the inverse operator. When the storage of many inverse operators (even as implicit matrices) is not feasible, a parameter-dependent preconditioner could be obtained by a projection into the linear span of  preconditioners, such as incomplete factorizations, sparse approximate inverses, H-matrices or other preconditoners that are readily available for a considered application. Also, we have restricted the presentation to the case of real matrices but the methodology can be naturally extended to the case of complex matrices.

\end{document}